\documentclass[reqno]{amsart}
\theoremstyle{plain}
\usepackage{latexsym}
\usepackage{amssymb}
\usepackage{epsfig}
\def\Xint#1{\mathchoice
   {\XXint\displaystyle\textstyle{#1}}%
   {\XXint\textstyle\scriptstyle{#1}}%
   {\XXint\scriptstyle\scriptscriptstyle{#1}}%
   {\XXint\scriptscriptstyle\scriptscriptstyle{#1}}%
   \!\int}

\def\XXint#1#2#3{{\setbox0=\hbox{$#1{#2#3}{\int}$}
     \vcenter{\hbox{$#2#3$}}\kern-.5\wd0}}

\def\dashint{\Xint-}

\newcommand{\chara}{1\!\!1}

\newcommand{\Zb}{\mathbb{Z}}

\newcommand{\lt}{\left}
\newcommand{\rt}{\right}

\newcommand{\nn}{\nonumber}

\newcommand{\lm}{\lambda}

\newcommand{\qd}{\quad}

\newcommand{\veps}{\varepsilon}

\newcommand{\vs}{\varsigma}

\newcommand{\vp}{\varphi}
\newcommand{\ep}{\epsilon}

\newcommand{\wt}{\widetilde}

\newcommand{\YI}{\mathcal{Y}}
\newcommand{\GI}{\mathcal{G}}

\newcommand{\SI}{\mathcal{S}}
\newcommand{\RI}{\mathcal{R}}

\newcommand{\DI}{\mathcal{D}}
\newcommand{\CI}{\mathcal{C}}
\newcommand{\HI}{\mathcal{H}}
\newcommand{\UI}{\mathcal{W}}
\newcommand{\OI}{\mathcal{O}}

\newcommand{\ti}{\tilde}

\newcommand{\BI}{\mathcal{B}}

\newcommand{\R}{\mathrm {I\!R}}

\newcommand{\ca}[1]{\mathrm{Card}\lt(#1\rt)}
\newcommand{\dia}{\diamondsuit}
\newcommand{\rest}{\llcorner}

\newcommand{\apep}{\ep}

\newtheorem{a1}{Lemma}
\newtheorem{a2}{Theorem}

\newtheorem{a5}{Proposition}
\newtheorem{a6}{Corollary}

\theoremstyle{remark}
\newtheorem{remark}{Remark}

\begin{document}
\title[On multiwell Liouville theorems in higher dimension] {On multiwell Liouville theorems in higher dimensions}
\author{Robert L. Jerrard \and Andrew Lorent}
\address{Department of Mathematics, University of Toronto,
Toronto, Ontario, Canada}\email{rjerrard@math.toronto.edu}
\address{Centro De Giorgi, Scuola Normale Superiore,
Piazza dei Cavalieri, Pisa, Italy}\email{andrew.lorent@sns.it}

\begin{abstract}
We consider certain subsets of the space of $n\times n$ matrices of
the form $K = \cup_{i=1}^m SO(n)A_i$, and we prove that
for $p>1, q \geq 1$ and for
connected $\Omega'\subset\subset\Omega\subset \R^n$, there exists
positive constant $a<1$ depending on $n,p,q, \Omega, \Omega'$ such
that for $ \veps=\| \mbox{dist}(Du, K)\|_{L^p(\Omega)}^p$ we have
$\inf_{R\in K}\|Du-R\|^p_{L^p(\Omega')}\leq M\veps^{1/p}$ provided
$u$ satisfies the inequality $\| D^2 u\|_{L^q(\Omega)}^q\leq
a\veps^{1-q}$.  Our main result holds whenever $m=2$, and also for
{\em generic} $m\le n$ in every dimension $n\ge 3$, as long as the wells
$SO(n)A_1,\ldots, SO(n)A_m$  satisfy a
certain connectivity condition. These conclusions are mostly known
when $n=2$, and they are new for $n\ge 3$.
\end{abstract}

\maketitle

\qd\qd\qd\qd\qd\qd\qd\qd\qd\qd\qd\qd\qd\qd\qd \today

\section{Introduction}

Rigidity theorems for mappings whose gradient lie in a subset of the conformal matrices date  back to 1850, when Liouville \cite{lou} proved that given a domain $\Omega\subset \R^n$ and a function
$u\in C^3\lt(\Omega,\R^3\rt)$ with the property that $Du\lt(x\rt)=\lambda\lt(x\rt) O\lt(x\rt)$
where $\lambda\lt(x\rt)\in \R_{+}$ and $O\lt(n\rt)\in SO\lt(n\rt)$ then $u$ is either affine or a Mobius
transformation.
A corollary to Liouville's Theorem is that a $C^3$ function whose gradient belongs everywhere to $SO\lt(n\rt)$ is an
affine mapping.
A striking quantitative version of this fact
was recently proved by
Friesecke, James and M\"uller \cite{fmul},
who showed that for every bounded open connected Lipschitz domain $U\subset \R^n$, $n\ge 2$,
and every $q>1$, there exists a constant $C(U,q)$
such that, writing $K := SO(n)$,
\begin{equation}
\inf_{R\in K} \label{mthm1} \| Dv-R\|_{L^q\lt(U\rt)}\leq C\lt(U,q\rt)\|
d\lt(Dv,K\rt)\|_{L^q\lt(U\rt)}
\quad\mbox{ for every } v\in W^{1,q}(U;\R^n).
\end{equation}
Here and below, $d(M,K)$ denotes the distance from a matrix $M\in \R^{n\times n}$
to a subset $K\subset \R^{n\times n}$, measured in the Euclidean norm.
This result strengthens earlier work of a series of authors, including
John \cite{fg1},\cite{fg1}, Reshetnyak \cite{res}, and Kohn
\cite{ko}, and
it has had a number of important applications. For example, it is a
main tool used to provide
a relatively complete analysis of the gamma limit of thin elastic structures, \cite{fmul}, \cite{fmul2}.

A number of works have extended the above result \eqref{mthm1}
to cover various larger classes of matrices than $SO\lt(n\rt)$. Faraco and Zhong proved
the corresponding result with
$K=\Pi SO\lt(n\rt)$ where $\Pi \subset \R_{+}\backslash \lt\{0\rt\}$ is a compact set,
\cite{fac}.
Chaudhuri and M\"uller \cite{chau} and later De Lellis and  Szekelyhidi \cite{las1} considered a set of the form $K=SO\lt(n\rt)A\cup SO\lt(n\rt)B$ where $A$ and $B$ are \em strongly
incompatible \rm in the sense of Matos \cite{matos}.

If we consider two compatible wells
$K =SO\lt(n\rt)A\cup SO\lt(n\rt)B$, i.e.\ wells for which there exists matrices $X\in SO\lt(n\rt)A$, $Y\in SO\lt(n\rt)B$ with
$\mathrm{rank}\lt(X-Y\rt)=1$, then the example of a piecewise affine function $u$
such that $\mathrm{Image}(Du) = \{ X,Y\}$ shows that no exact analog of \eqref{mthm1} can hold. In this paper we show, however, that a sort of $2$-well theorem
can hold provided one has suitable control over second derivatives; indeed, this remains
true for collections of $m\ge 3$ wells $K = \cup_{i=1}^m SO(n)A_i$  satisfying certain algebraic conditions. As we will recall in greater detail below, most of our main conclusions are known in $2$ dimensions, however all are new in $\R^n, n\ge 3$. The main result of this paper is

\begin{a2}

\label{T1}
Let $p,q\geq 1$, let $A_1,A_2,\dots A_m\in \R^{n\times n}$ be matrices of non-zero determinant,
and let $K=\bigcup_{i=1}^m SO\lt(n\rt)A_i$.
Suppose that
$m=2$,
or that for each $i\in \{1,\ldots, m\}$, there exists $v_i\in S^{n-1}$ such that
either
\begin{equation}
|A_i v_i| > |A_j v_i|\mbox{ for all }j\ne i
\label{T2.h1}
\end{equation}
or
\begin{equation}
|v_i^TA_i^{-1} | > |v_j^T A^{-1}_j|\mbox{ for all }j\ne i.
\label{T2.h2}
\end{equation}
Then for any bounded, open, $\Omega\subset \R^n$
and connected $\Omega'\subset\subset\Omega$
there exists positive constants $a<1$ and $M<\infty$, depending on
$K, \Omega, \Omega', p,q$,  such that for any $u\in W^{1,p}\cap W^{2,q}(\Omega;\R^n)$
that satisfies
\begin{equation}
\label{ee3} \frac 1\vs \int_{\Omega} d^p\lt(Du,K\rt) + \vs^q\lt|D^2
u\rt|^q dx \leq  a
\end{equation}
for some $\vs\in (0,1]$,
there exists $i\in\lt\{1,\dots m\rt\}$ such that
\begin{equation}
\label{xx670}
\int_{\Omega'} d^p\lt(Du,SO\lt(n\rt)A_i\rt) dx\leq M\vs^{1/p},
\end{equation}
and if $p>1$ there exists $R\in SO\lt(n\rt) A_i$ such that
\begin{equation}
\label{wem1.8}
\int_{\Omega'} \lt|Du-R\rt|^p dx\leq M \vs^{1/p}.
\end{equation}
\end{a2}

The theorem is interesting in when $0< \vs \ll a$. The result as stated follows easily from
the case when $\Omega$ is the unit ball in $\R^n$ and $\Omega'$ is some small subball,
so we will mostly focus on this situation.
The
conclusions of the theorem are generally not true if $\Omega' = \Omega$, as long as compatible wells are allowed; this is easily seen by taking $u$ to be a suitable mollification of a piecewise affine function whose gradient assumes exactly two values. An example in \cite{contswb}, Remark 6.1, shows that the scaling in \eqref{xx670}, \eqref{wem1.8} is sharp.

\begin{remark}
We suspect that the theorem  remains true whenever $m=3, n\ge 2$, and
we verify in Section \ref{sctwt} that for $m = n \ge 3$, the hypotheses of the theorem
are {\em generically} satisfied as long as the $n$ wells have the property that they cannot be
partitioned into two disjoint subfamilies of wells with no rank-$1$ connections between them.

However, for $m=4$  and any $n\ge 2$,
one can find examples of matrices $A_1,\ldots, A_4$ such that
the conclusions of the theorem fail for $K = \cup_{i=1}^4 SO(n) A_i$.
To construct an example for  $\Omega' \subset\subset \Omega \subset \R^2$,
we start with a equilateral triangle $T\subset \Omega'$ of diameter $\ell$, and we partition $T$ into three congruent subtriangles $S_1,S_2, S_3$. Let $S_4 = \Omega\setminus T$.
We can then find a piecewise affine function $u_0$ and matrices
$A_1,\ldots A_4\in M^{2\times 2}$
such that $Du_0 = A_i$ a.e. on $S_i$, for $i=1,\ldots,4$.
Let $u_\vs =u_0*\phi_{\ep}$ where
$\phi_{\ep}:=\ep^{-n} \phi\lt(\frac{z}{\ep}\rt)$ and $\phi$ is a
standard mollifier on $\R^2$.

One can fix $\ell \lessapprox a$  such that $u_\vs$
satisfies \eqref{ee3} for every $\vs \ll l$. However, as $\vs\to 0$,
$\int_{\Omega'} d^p(u_\vs, SO(n)A_i)dx  \gtrapprox c \ell^2$
for every $i$, so the conclusions of the theorem do not hold.


\end{remark}

\vspace{.1in} The first 2-well Liouville Theorem was due to the
second author \cite{lor2}, who established essentially the above
result in the case when $m=n=2$ and $p=q=1$, for matrices $A,B\in
\R^{2\times 2}$ with $\det A =\det  B$, with suboptimal scaling in
\eqref{wem1.8}, and under the  assumption that $u$ is Lipschitz and
invertible, with Lipschitz inverse. This was greatly improved by
Conti and Schweizer, \cite{contsw}, who proved Theorem \ref{T1} for
$q=1$, still for $m=n=2$. In particular \cite{contsw} established
this case of the theorem with the optimal scaling as in
\eqref{xx670}, \eqref{wem1.8}, and without either the assumption of
invertibility or any conditions on the two wells. A different proof
of Theorem \ref{T1}  for $n=m=2$, valid  for general $p,q\ge 1$, was
given in \cite{lor6}. This argument is conceptually simple, and the
proof clarifies some technical issues in \cite{contsw}, but it
yields suboptimal scaling in \eqref{wem1.8} and requires the
assumption $\det A = \det B$.



\subsection{Ingredients in the proof}

As mentioned above, we work mostly on $\Omega = B_1\subset \R^n$.
Straightforward arguments from previous work, recalled in Section \ref{S2}, allow us easily to find some $i_*\in \{1,\ldots, m\}$
and a large set $U_0\subset B_1$ with small perimeter in $B_1$, such that
$d(Du, SO(n)A_{i_*}) = d(Du,K)$ in $U_0$.
We always assume for concreteness that $i_*=1$.
Our first goal is to find many pairs of points $(x,y)\subset U_0\times U_0$
such that
\begin{equation}
|u(x) - u(y)| = |A_1(x-y)| + O(\vs^{1/p}).
\label{firstgoal}\end{equation}
Further easy and well-known arguments, also recalled in Section \ref{S2}, allow us to find many pairs of
points for which the inequality $|u(x) - u(y)| \le |A_1(x-y)| + O(\vs^{1/p})$ holds. Following previous work,
we wish to prove the opposite inequality by applying the same argument to $u^{-1}$. In general $u$ is not
invertible, but in fact it is only necessary to prove that there are many line segments along which $u$ can be inverted.
One of the important contributions of \cite{contsw} was to introduce arguments, using tools from degree theory, to support this
contention.
Their local invertibility arguments, however, rest on the Sobolev embedding $W^{2,1}\hookrightarrow H^1$
(in ways that are not made completely explicit), and so do not apply to $\R^n$ for $n\ge 3$.

To address this difficulty we prove a new Lipschitz truncation result, showing that one can find a
Lipschitz function $w$ such that the set $\{ x\in B_1 : w(x)\ne u(x)\}$ is not only small,
but also can be contained in a set of small perimeter. The new point is the perimeter estimate, which follows
from the control over second derivatives of $u$ supplied by \eqref{ee3}.
The specific facts we need about this Lipschitz approximation are proved in Section \ref{S:lip}. They are deduced from a
general  truncation result that we  prove  in Section \ref{S:app}.
Using the Lipschitz approximation and some elements from earlier work of various authors,
we find in Section \ref{S4} a large subset $\UI$ of $u(B_1)$
on which an inverse map is well-defined and Lipschitz, with its gradient near $A_1^{-1}SO(n)$ and, crucially,
with control over the perimeter of $\UI$. This allows us in
Section \ref{S5} to complete the proof that \eqref{firstgoal} holds for a  large set of pairs of points.

The proof of Theorem \ref{T1} is given in Section \ref{S6}. We first consider the case when
the majority phase, represented by $A_1$, satisfies \eqref{T2.h1}.
Then we can bound $d(\cdot, SO(n)A_1)$ by $d(\cdot,K) + \mbox{a null lagrangian}$,  and it directly follows, via integration by parts, that
\begin{equation}
\sum_{k=0}^n \int_{[x_k,x]} d(Du , SO(n)A_1) dH^1
\le  C \sum_{k=0}^n\int_{[x_k,x]} d(Du , K) dH^1 +  \mbox{ boundary terms}
\label{2ndgoal}\end{equation}
where $x_0,\ldots, x_n$ are the vertices of a long, thin simplex with long axis
roughly parallel to $v_i$,  $x$ is a point near the barycenter, and $[x_k,x]$ denotes
the line segment joining $x_k$ and $x$. The boundary terms have the form
$C\sum_{k=0}^n |u(x_k)- l_R(x_k)|$, where $l_R$ is an affine map with
$Dl_R = R \in SO(n)A_1$.
The inequality \eqref{2ndgoal} recasts and extends
ideas developed in \cite{contsw} for $n=2$.
We present the short proof of \eqref{2ndgoal} in the next subsection.

If the majority phase $A_1$ satisfies \eqref{T2.h1}, then by using
\eqref{firstgoal} and a linear algebra lemma
proved in Section \ref{S:app}, we can find a vertices $x_0,\ldots, x_n$
and an affine $l_R$ map such that the boundary terms in \eqref{2ndgoal} are less than $C\vs^{1/p}$. The proof of Theorem
\ref{T1} in this case is essentially completed by
integrating \eqref{2ndgoal} over points $x$ near the barycenter.

When the majority phase $A_1$ satisfies \eqref{T2.h2}, the idea of the proof is
to apply to $u^{-1}$  the argument already used to
prove the theorem under assumption \eqref{T2.h1}. The fact that $u$ need not
be invertible again causes technical difficulties.
Thus, we work with the Lipschitz
approximant $w$ found earlier, and we use a lemma, proved in Section \ref{S:app},
which asserts roughly speaking that {\em almost every}
line segment passing through a large convex subset of $w(B_1)$
can be realized as the image via $w$ of a Lipschitz path in $B_1$.
Although the restriction of $w$ to these Lipschitz paths is not injective
in general, this lemma provides a good enough proxy for invertibility to allow us
to complete the proof of the theorem under the hypothesis \eqref{T2.h2}.
The null lagrangian calculation that leads to \eqref{2ndgoal} is a bit harder
to implement in the inverse direction, and in its place we  use an argument more directly
related to a proof given in \cite{contsw} when $n=2$.

Finally,  it is easy to see that when $m=2$, each well must satisfy at least one of
\eqref{T2.h1}, \eqref{T2.h2}.

The condition \eqref{T2.h2} does not appear in any previous work, so that our result
yields new information when $m\ge 3$, even in $2$ dimensions.
In particular, in $2$ dimensions \cite{contsw} essentially proves the theorem if every
$A_i$ satisfies either
\eqref{T2.h1} or the condition that
\begin{equation}
\mbox{ for each $j\ne i$,
$\det A_j > \det A_i>0$}.
\label{T2.h3}\end{equation}
Only  the case of $m=2$ wells is discussed in \cite{contsw}, but the argument works
almost without change for $m\ge 2$ under the assumptions discussed here. The proof given under condition
\eqref{T2.h3} is intrinsically $2$-dimensional and so is not available here in the generality we consider here.

\subsection{Proof of \eqref{2ndgoal}}

As discussed above, a crucial point in the proof of Theorem \ref{T1}
in the case when hypothesis \eqref{T2.h1} holds
is that if $\{ x_0,\ldots , x_n\}$ are the vertices of a suitable
simplex (where ``suitability'' is related to the algebraic condition
\eqref{T2.h1})
then one can bound
$\int d(Du,SO(n)A_1)$ by $\int d(Du,K)$ +  {\em boundary terms}
along certain lines.
We illustrate how this works in the simplest possible case,
that of a $2$-well Liouville Theorem in $\R^1$.
For this, suppose that $K = \{ a_1, a_2 \}$ for $a_2<a_1\in \R$, and
consider  $u: (-1,1)\to \R$. Since $a_1>a_2$,  we can find constants $c_1, c_2$ such that
\begin{equation}
|s-a_1| = d(s, a_1) \le  c_1d(s,K) + c_2 (a_1- s)\quad\quad\mbox{ for all }s\in \R.
\label{easycase1}\end{equation}
We substitute $s=u'$ in \eqref{easycase1} and integrate. If we let $l_1$
be an affine function with $l_1' = a_1$, then
$a_1-u' = (l_1 - u)'$, and we find that
\begin{align}
\int_{-1}^1 d(u', a_1)
\le c_1\int_{-1}^1 d(u' ,K)   + c_2 \left( | l_1(1)- u(1)| + | l_1(-1) - u(-1)| \right).
\label{easycase2}\end{align}

The next lemma, which is not used until Section \ref{S6},  is essentially  the same argument, but now for
$m$ wells in $\R^n$.
Note that if $i=1$ satisfies \eqref{T2.h1}, then condition \eqref{shrink.dir} below
is fulfilled if $\{x_0,\ldots, x_n\}$ are the vertices of a long thin simplex roughly parallel to
$v_1$.

\begin{a1}
Assume
that $\{A_1,\ldots, A_m\}$ are $n\times n$ matrices 
and let $K = \bigcup_i SO(n) A_i$.
Let $x_0,\ldots, x_n\in B_1\subset \R^n$ be vertices of a simplex with
the property that
\begin{equation}
\lt| A_{1} \frac {x-x_i}{|x-x_i|} \rt| > (1+\alpha) \lt| A_j \frac {x-x_i}{|x-x_i|} \rt| \quad \mbox{ for all }j\in \{2,\ldots, m\}
\mbox{ and }i\in \{0,\ldots, n\}
\label{shrink.dir}\end{equation}
for some $x$ in the interior of the simplex
${\rm conv}\{ x_0,\ldots, x_n\}$. Then there exists a constant $C$ such that
\begin{equation}
\sum_{i=0}^n \int_{[x_i, x]} d(Du, SO(n) A_1) \, dH^1 \ \le  \ C
 \sum_{i=0}^n \int_{[x_i, x]} d(Du, K) \, dH^1 \ +
C \sum_{i=0}^n  |u(x_i) - l_R(x_i)|
\label{IbyP}\end{equation}
for every  smooth $u:B_1\to \R^n$ and every affine map $l_R$ with $Dl_R = R\in SO(n)A_1$.

Moreover, if we write $x = \sum_{i=0}^n \lambda_i x_i$ with $\sum
\lambda_i =1$ and $\lambda_i>0$ for all $i$, and if $\lambda_i
|x-x_i| \ge\alpha'>0$ for all $i$, then the constant $C$ in
\eqref{IbyP} are uniformly bounded by constants depending only on
$\{ A_i\}$, $\alpha$, $\alpha'$. \label{L.IbyP}\end{a1}

This lemma is inspired by an argument from
\cite{contsw}. In the context of the $1$-dimensional toy problem discussed above,
the idea in \cite{contsw} would be to use information about
$\int_{-1}^1 d(u',K)$ and the boundary behavior of $u$
at $\pm 1$ to bound $L^1( \{ x\in (-1,1) :
d(u', a_2) <d(u', a_1) \})$. We use arguments of this sort in Lemma \ref{LL2},
when considering the hypothesis \eqref{T2.h2}. In fact, either argument ---integration by parts or direct estimates of the size of a bad set ---  could be used
to prove both halves of Theorem \ref{T1}.

\begin{proof}[Proof Lemma \ref{L.IbyP}]
{\em Step 1}.
Fix $\{ x_0,\ldots, x_n\}$ and $x = \sum_{i=0}^n \lambda_i x_i$ satisfying \eqref{shrink.dir},
where $0 < \lambda_i$ for all $i$, and
$\sum \lambda_i=1$. Also fix an affine map $l_R$ with $Dl_R = R\in SO(n)A_1$.

For $i=0,\ldots, n$, let us write
$\tau_i := \frac {x-x_i}{|x-x_i|}$ and $v_i := \lambda_i R(x-x_i) = \lambda_i |x-x_i|  R\tau_i$.
Note that $\sum v_i = R( \sum \lambda_i(x-x_i)) = 0$.

We first claim that \eqref{shrink.dir} implies that there exist $c_1, c_2>0$ such that
\begin{equation}
d( M, SO(n)A_1) \ \le  c_1 d(M , K) +  c_2 v_i^T(R-M)\tau_i
\label{hardcase1}\end{equation}
for every $n\times n$ matrix $M$ and
every $i\in \{0,\ldots,n\}$.
Here  $v_i,\tau_i$ are  column vectors, and $v_i^T$ denotes the transpose of $v_i$.
Inequality (\ref{hardcase1}) is the analog of \eqref{easycase1} from the $1$-dimensional case.
To prove \eqref{hardcase1}, we write $\tilde \lambda_i = \lambda_i|x-x_i|$ for simplicity, and we
note that since $R\in O(n)A_1$,
\[
v_i^T R \tau_i = \tilde \lambda_i  |R\tau_i|^2
= \tilde \lambda_i |A_1\tau_i|^2
> \tilde \lambda_i(1+\alpha) |A_1\tau_i| \  |A_j\tau_i|
\]
for $j\ge 2$, using \eqref{shrink.dir}.
Similarly
$v_i^T M \tau_i \le \tilde \lambda_i |A_1\tau_i| \ |M \tau_i|$. In
particular,
if $M\in SO(n) A_j$ for some $j\ge 2$, then $|M \tau_i| = |A_j\tau_i|$,
and so
\[
v_i^T (R- M) \tau_i \ge \tilde \lambda_j \alpha |A_1\tau_i| \ |A_j \tau_i| \ge c>0 \ \quad
\mbox{ for } \ \ M\in \cup_{j=2}^m SO(n)A_j.
\]
Also, if $M\in \cup_{j\ge 2}SO(n)A_j$, then $d(M, SO(n) A_1) \le C(\{A_1,\ldots, A_m\})$,
It follows that we can fix positive constants $c_2$ so large and $\delta$ so small that
\[
d( M, SO(n)A_1) \ \le  c_2 v_i^T(R-M)\tau_i
\]
say for all $M$ such that $d(M,\cup_{j=2}^m SO(n) A_j) \le \delta$.
Then by choosing $c_1$ large enough, we can arrange that
\[
d( M, SO(n)A_1)  -   c_2 v_i^T(R-M)\tau_i
\ \le  c_1 d(M, K)
\]
whenever $d(M,\cup_{j=2}^m SO(n) A_j) \ge \delta$. Then \eqref{hardcase1} follows.

{\em Step 2}. Now we substitute $M= Du$ in \eqref{hardcase1}, so that $R-M$ becomes $R-Du = D(l_R-u)$. We then integrate
to find that
\begin{equation}
\sum_{i=0}^n \int_{[x_i,x]}
d( Du, SO(n)A_1) \, dH^1
\le \    \sum_{i=0}^n \int_{[x_i,x]}
\left[ c_1 d( Du, K) +  c_2 v_i^T  (R -  D u)\tau_i\right] \ dH^1
\label{hardcase2.1}\end{equation}
Since $R-Du = D(l_R-u)$ and $\tau_i$ is tangent to $[x_i,x]$, we can integrate
by parts to find that
\begin{align}
\sum_{i=0}^n \int_{[x_i,x]}
v_i^T  (R -  D u)\tau_i \ dH^1
&=
\sum_i v_i^T\left[(l_R - u)(x) - (l_R - u)(x_i)\right]\nn \\
&=
- \sum_i v_i^T (l_R - u)(x_i)
\label{hardcase2.2}
\end{align}
since $\sum_i v_i^T(l_R - u)(x) = \left( \sum_i v_i  \right)^T (l_R-u)(x) = 0$.
Now \eqref{IbyP} follows by combining \eqref{hardcase2.1},
\eqref{hardcase2.2}.

The statement about dependence of the constants in \eqref{IbyP}
on various other parameters
follows from inspection of the above argument.
\end{proof}

\section{Preliminaries}
\label{S2}

In this section we introduce some notation and reformulate some arguments
from \cite{lor2} that provide the starting point for our analysis.

\subsection{Some notation}
Given matrices $A_1,\ldots, A_m$, we always write $K = \cup_{i=1}^m SO(n)A_i$.

We write $B_r(x)$ for the {\em open} ball in $\R^n$ of radius $r$, centered at $x$. We write $B_r$ as an
abbreviation  for $B_r(0)$. Define $[x,y]$ to denote the line segment joining $x$ and $y$.
If $S$ is a subset of $\R^n$, then $\chara_S$ always denotes the characteristic function of $S$,
so that $\chara_S(x) = 1$ if $x\in S$ and $0$ otherwise.

We will write $\sigma= \sigma(K)$ to denote a fixed small number depending only
on the given matrices $A_1,\ldots, A_m$. We select $\sigma\le 1$ to satisfy
\begin{equation}
\sigma < \frac 14 dist(SO(n)A_i, SO(n)A_j)\quad\mbox{ for all }i\ne j,
\label{sigma.c1}\end{equation}
\begin{equation}
d(M, K) > \sigma \quad\mbox{ for any matrix $M$ such that }\det M < \sigma.
\label{sigma.c2}\end{equation}

Note that \eqref{sigma.c1} implies that
$d(Du, K) = d(Du,  SO(n)A_i)$ whenever $d(Du, SO(n)A_i)<2\sigma$.

All constants throughout, including generic constants $C$ that appear in many estimates, as well as
named constants such as $\kappa_0$ in Proposition \ref{L9} for example, may depend on the collection
$K$ of wells,  the dimension $n$, and the powers $p,q$ appearing in assumptions \eqref{ee3}, for example.
but are independent of the parameters $\vs, a$.

We often (though not always) use latin letters to refer to the reference configuration $B_1$
and greek letters to refer to the image $u(B_1)$.
Thus points in $B_1$ will be denoted $x,y,z$, whereas points in the image will be
denoted $\xi,\eta, \zeta$. In addition we will write $\beta_\rho$ to denote an ellipsoid in the image
with length-scale $\rho$; in fact $\beta_\rho$ will be defined as
$\beta_\rho = l_R(B_\rho)$, where $l_R$ is a particular
affine map we find that is close to $u$, see Section \ref{S:lip}.

\subsection{Finding a majority phase}

\begin{a1}
\label{L3}
Let $K=\cup_{i=1}^m SO\lt(n\rt)A_i $.
Let $u:B_1\subset \R^n\to \R^n$ be a smooth function such that
\begin{equation}
\label{e1}
\frac 1 \vs \int_{B_1}\left(  d^p\lt(Du,K\rt) + \vs^q \lt|D^2 u\rt|^q \right) dx \le a.
\end{equation}
Then we can find $i\in\lt\{1,\ldots, m\rt\}$ an open set $U_0\subset B_1$
with smooth boundary
such that
\begin{equation}\label{L3.c1}
\mathrm{Per}_{B_1}\lt(U_0\rt)<C a\text{ and }L^n\lt(B_1\backslash
U_0\rt)<C a^{\frac{n}{n-1}},
\end{equation}
and 
\begin{equation}\label{a42}
d(Du,SO(n)A_i) = d(Du,K) <\sigma
\quad \quad
\mbox{ for all }x\in U_0.
\end{equation}
\end{a1}

We take $U_0$ to be smooth because it is convenient later to identify
$\mathrm{Per}_{B_1}\lt(U_0\rt)$  with $H^{n-1}(B_1\cap \partial U_0)$.

\begin{proof}[ Proof of Lemma \ref{L3}. ]
Let $q^*$ be the Holder conjugate of $q$, and let
$s:= 1+\frac{p}{q^*}$ and
$J\lt(x\rt):=d^{`s} \lt(Du\lt(x\rt),K\rt)$.
If $q^*=\infty$
we use the convention that $d^{\frac{p}{q^*}}\lt(\cdot,K\rt):=1$. We have by Young's inequality
\begin{equation}
\int_{B_1}|DJ| \ dx  \le C \int_{B_1}  d^{\frac{p}{q^*}}\lt(Du,K\rt)\lt|D^2 u\rt| dx
\ \le \  \frac C\vs \int_{B_1} d^p\lt(Du,K\rt)+\vs^q\lt|D^2 u\rt|^q dx
\overset{\eqref{e1}}{\le} \  Ca.
\label{youngsineq}\end{equation}
Then by the coarea formula, we can find
$\alpha\in \lt( (\frac {\sigma}2)^{s},
{\sigma^{s} }
\rt)$ with
$Per_{B_1}( \{ x \in B_1: J(x)<\alpha\}) \le C a$.
Note that
\[
\bigcup_{i=1}^n  \{x\in B_1:d^s\lt(Du,SO\lt(n\rt)A_i\rt) < \alpha\}
=
\{ x\in B_1 : J(x)<\alpha\}.
\]
Since the sets on the left-hand side above are disjoint by the choice \eqref{sigma.c1} of $\sigma$, it follows that
\[
Per_{B_1} (\lt\{x\in B_1:d^s\lt(Du,SO\lt(n\rt)A_i\rt) < \alpha\rt\} )  \ \le  \ \
Per_{B_1}( \{ x \in B_1: J(x)<\alpha\}) \ \ \le \ \  C a
\]
for every $i$.
So by the relative isoperimetric inequality we have
\begin{align*}
&\min\lt\{L^n\lt(\lt\{x\in B_1:d^s\lt(Du,SO\lt(n\rt)A_i\rt)>\alpha\rt\}\rt),
L^n\lt(\lt\{x\in B_1:d^s\lt(Du,SO\lt(n\rt)A_i\rt)<\alpha\rt\}\rt)\rt\}\\
&\quad\quad\qd\qd\qd\qd\qd<C a^{\frac{n}{n-1}}
\end{align*}
for every $i$. Since $\int d^p(Du, K)< a \vs$, it cannot be the case that
$\lt\{d^s\lt(Du,SO\lt(n\rt)A_i\rt)<\alpha\rt\}$ has small measure for
every $i$, and since $a$ is small, there can be at most one $i$ such that
$L^n\lt(\lt\{x\in B_1:d^s\lt(Du,SO\lt(n\rt)A_i\rt)<\alpha\rt\}\rt) > 1 - C a^{\frac n{n-1}}$.
We define
$$
U_0 :=\lt\{x\in B_1:d^s\lt(Du,SO\lt(n\rt)A_i\rt)<\alpha\rt\}\text{ for this choice of }i.
$$
Since $J$ is a $C^1$ function by Sard's Theorem the image under $J$
of the critical points of $J$ have zero $L^1$ measure, so we can
assume we choose $\alpha$ so that the level set
$J^{-1}\lt(\alpha\rt)$ does not intersect the set of critical points
of $J$. Then $\partial U_0$ is smooth, as required.

\end{proof}

Upon relabeling, we may assume that $i=1$
in Lemma \ref{L3}, so that  $U_0$ satisfies
\begin{equation}
\label{e3}
U_0 \subset \lt\{x\in B_1:d\lt(Du\lt(x\rt),SO \lt(n\rt)  A_1  \rt)<\sigma\rt\}.
\end{equation}
It would of course be possible to perform a change of variables that sets $A_1$ equal to the identity matrix.
We will mostly remain in the original coordinates, so that one can see explicitly where  $A_1$ appears in our arguments.

\subsection{Non-stretching pairs}

We next show that can find many pairs of points that are not stretched by $u$ (relative to the affine maps with
gradient in $SO(n)A_1$).
The argument we give is somewhat more complicated than necessary for the present lemma,
but it will be needed again in Section \ref{S5}.

\begin{a1}
\label{L6} Assume $u:B_1\subset \R^n\to \R^n$ is a smooth function that satisfies (\ref{e1})
and that $A_1$ is the majority phase as in (\ref{e3}).
Then
there exists $\GI_1\subset B_1\times B_1$
such that
\begin{equation}
L^{2n}( (B_1\times B_1)\setminus \GI_1) \le C a^{1/p}
\label{G1.1}\end{equation}
and letting $\ep=\vs^{\frac{1}{p}}$,
\begin{equation}
\mbox{ if }(x,y)\in \GI_1, \mbox{ then }
|u\lt(y\rt)-u\lt(x\rt)|\leq |A_1(y-x)|+ C \ep.
\label{G1.2}\end{equation}
\end{a1}

\begin{proof}
We define
\begin{equation}
\GI_1 :=
\{ (x,y)\in  B_1\times  B_1\ : \ [x,y] \subset U_0, \ \  \int_{[x,y]} d(Du, SO(n) A_1) \ dH^1 \ \le \ep\}.
\label{G1.def}
\end{equation}
Note from (\ref{e1}) we have
\begin{equation}
\label{wieb1}
\int_{B_1} d\lt(Du,K\rt)
\leq
C\lt(\int_{B_1} d^p\lt(Du,K\rt)\rt)^{\frac{1}{p}}
\leq
C a^{\frac{1}{p}}\ep.
\end{equation}

{\em Step 1} .
To prove \eqref{G1.2}, we fix $(x,y)\in \GI_1$, and we write
$\tau := \frac{y-x}{|y-x|}$.
Note that if $M$ is any $n\times n$ matrix, then,
$|M\tau| \le |A_1\tau| + C d(M, SO(n) A_1 )$.
Thus
\begin{align*}
|u(y) - u(x)|
&=
\left| \int_{[x,y]} Du(z) \ \tau \  dH^1z \right| \\
&\le
\int_{[x,y]} \left[|A_1\tau|+ C \,d(Du,SO(n)A_1) \right] dH^1
\ \le \
|A_1(x-y)| + C \ep
\end{align*}
for $(x,y)\in \GI_1$.

{\em Step 2}.
We next prove \eqref{G1.1}.
If $(x,y)\in  \left(B_1 \times B_1 \right)\setminus \GI_1$,
then at least one of the following must hold: either $x$ or $y$ fails to belong to $U_0$, that is
\begin{equation}
(x,y) \in  [B_1\setminus U_0] \times B_1\text{ or } (x,y) \in  B_1\times [B_1\setminus U_0];
\label{notG11}\end{equation}
or the segment $[x,y]$ meets $\partial U_0\cap B_1$, that is
\begin{equation}
[x,y]\cap ( \partial U_0 \cap B_1 )\ \ne \  \emptyset;
\label{notG12}\end{equation}
or
\begin{equation}
\int_{[x,y]}  \chara_{U_0} \ d(Du, SO(n)A_1) \ dH^1 \ > \ep.
\label{notG13}\end{equation}
We saw in Lemma \ref{L3} that $L^n(B_1\setminus U_0)\le Ca^{\frac n{n-1}}$,
so clearly \eqref{notG11} holds on a set of measure at most
$Ca^{\frac n{n-1}}\le Ca$.
And Lemma \ref{L7.7} (proved at the end of this section)
shows that
\[
L^{2n}(\{(x,y) \ : \ \mbox{ \eqref{notG12} holds } \} ) \le C
H^{n-1} ( \partial  U_0\cap B_1 ).
\]
However, in Lemma \ref{L3} we showed that
$H^{n-1} (  \partial  U_0\cap B_1  ) \le Ca$.
Finally,  Lemma \ref{L.goodpairs} (proved immediately below) implies that
\begin{align*}
L^{2n}(\{(x,y)\in B_1\times B_1 \ : \ \mbox{ \eqref{notG13} holds } \} )
\ \ &\le\ \
\frac{C}{\ep} \int_{B_1}  \chara_{U_0} \ d(Du, SO(n)A_1) \\
\ \ &\overset{\eqref{a42}} {\le} \ \
\frac {C}{\ep}\int_{B_1} d(Du,K) \
\ \ \overset{\eqref{wieb1}}{\le} \ \  Ca^{1/p}.
\end{align*}
Together, these estimates imply \eqref{G1.1}.

\end{proof}

The first of the lemmas used above is

\begin{a1}
Suppose that $f:B_1\to \R$ is nonnegative and integrable. Then for any constant $b>0$,
\[
L^{2n} \left( \{ (x,y)\in B_1\times B_1 : \int_{[x,y]} f \ dH^1 > b \} \right) \ \le \ \frac C b \int_{B_1}f.
\]
\label{L.goodpairs}\end{a1}

\begin{proof}
We extend $f$ by $0$ on the complement of $B_1$. Then
by  a change of variables, we find that
\begin{align*}
\int_{B_1} \int_{B_1} \int_{[x,y]} f  \ dH^1 \ dy\  dx
&=
\int_{B_1}\int_{B_1}\int_0^1 f(x+ s(y-x)) |y-x| \ ds\  dy\  dx\\
&\le
\int_{B_1}\int_{B_2}\int_0^1 f(x+ sp) |p| \ ds\  dp\  dx.
\end{align*}
But
$\int_{B_1}\int_{B_2}\int_0^1 f(x+ sp) |p| \ ds\  dp\  dx \le
\| f \|_{L^1(B_1)}\int_{B_2}\int_0^1  |p| \ ds\  dp = C \|f\|_{L^1(B_1)}$ by Fubini's Theorem, so
the lemma follows by Chebyshev's inequality.
\end{proof}

The other lemma we used is

\begin{a1}
\label{L7.7}
Suppose that $\Omega$ is a bounded, open convex subset of $\R^n$.
Then there exists  a constant $C=C(\Omega)$ such that
for any set $S\subset \R^n$,
\begin{equation}
L^{2n}\left(
\{ (x,y)\in \Omega\times \Omega : [x,y]\cap S \ne \emptyset
\right) \le C H^{n-1}_\infty(S)
\label{goodpairs2}\end{equation}
where
$H^{n-1}_\infty\lt( S \rt) := \inf\{ \sum_i\gamma_{n-1} s_i^{n-1}  : \
S \subset \bigcup_{i} B_{s_i}\lt(x_i\rt)\}$.
\end{a1}

The constant $\gamma_{n-1}$ appearing in the definition of $H^{n-1}_\infty$ is the same normalization factor
appearing in the definition of Hausdorff measure, so that
$H^{n-1}_\infty(S)\le H^{n-1}(S)$ for every $S$.

\begin{proof}[ Proof of Lemma \ref{L7.7} ]

Without loss of generality we assume $0\in\Omega$. For any $S\subset \R^n$, we will write
\[
\vp(S) := L^{2n}\left(
\{ (x,y)\in \Omega\times \Omega : [x,y]\cap S \ne \emptyset \}
\right).
\]
We first claim that
$\vp(B_r(p)) \le C r^{n-1}$ for any $p\in \R^n$ and $r>0$.
To prove this, note that by Fubini's Theorem,
\begin{align*}
\vp( B_r(p))
&= \int_{y\in \Omega} L^n( \{ x\in \Omega \ : [x, y]\cap B_r(p)\ne \emptyset\} \ )\ dy \\
&\le  \int_{|z|\le \mathrm{diam}(\Omega)} L^n( \{ x\in \Omega \ : [x, x+z]\cap B_r(p)\ne \emptyset\} \ ) \  dz.
\end{align*}
And for every fixed $z\ne 0$, if $[x,x+z] \cap B_r(p) \ne \emptyset$, then $x$ belongs to
the cylinder of radius $r$, with axis parallel to $z$, that contains $B_r(p)$. The intersection of such
a cylinder with $\Omega$ has $L^n$ measure bounded by $Cr^{n-1}$. Thus
$L^n( \{ x\in \Omega \ : [x, x+z]\cap B_r(p)\ne \emptyset\} \ ) \le Cr^{n-1}$ for every $z\ne 0$. The claim follows.

Now given $S\subset \R^n$, let $\{ B_{r_i}(p_i) \}$ be a collection of balls such that
\[
S\subset \bigcup_i B_{r_i}(p_i) \quad\quad\mbox{ and }\quad
\sum_{i} \gamma_{n-1} r_i^{n-1} \le 2 H^{n-1}_\infty(S).
\]
Then any segment $[x,y]$ that intersects $S$ also intersects some ball $ B_{r_i}(p_i)$, so we deduce that
\[
\vp( S )
\le \sum_i \vp(B_{r_i}(p_i) ) \le C \sum_i r_i^{n-1} \le C H^{n-1}_\infty(S).
\]
\end{proof}

%
%

\section{Lipschitz approximation}\label{S:lip}

In this section we find a Lipschitz function $w$
that agrees with $u$ on the complement of
a small set $E$ and that is close to affine if $a$ is small. Such arguments are standard.
The main new ingredient here, which is crucial for our later arguments, is that we
use information about the second derivatives of $u$ to control
the perimeter of the set $E = \{ x: u(x) \ne w(x)\}$, or more precisely,
of a set that contains $E$.

\begin{a5}
\label{L9} Suppose the smooth function $u:B_1\subset \R^n\to \R^n$ satisfies (\ref{e1}),
and assume as in \eqref{e3} that $SO(n)A_1$ is the majority phase.

Then there is a Lipschitz function
$w:B_1\to \R^n$ with $\|Dw\|_{L^{\infty}\lt(B_1\rt)}\leq
C(K)$ and an open set $U_1\subset B_1$ with countably piecewise smooth boundary such that
$u=w$ in $U_1$,
and letting $\ep=\vs^{\frac{1}{p}}$ the following hold:

\begin{align*}
\it (i) \rm\quad\quad  &\|d\lt(Dw,K\rt)\|_{L^1(B_1)} \leq C \apep.\\
\it (ii) \rm\quad\quad
&d\lt(Dw\lt(x\rt),SO\lt(n\rt)A_1 \rt)=d\lt(Dw\lt(x\rt),K\rt)\text{ for
every }x\in U_1.\\
\it (iii) \rm\quad\quad &L^n\lt(B_1\backslash U_1\rt)\leq c
a^{\frac{n}{n-1}},\text{ and }Per_{B_1}\lt(U_1\rt)\leq c a.
\\
\it (iv) \rm\quad\quad &\text{There exists }R\in SO\lt(n\rt)A_1 \text{ and an
affine map }l_R\text{ with }Dl_R=R\text{ such that }\\
&\quad\quad\quad \|w-l_R\|_{L^{\infty}\lt(B_1\rt)}\leq c
a^{\frac{1}{n+1}}.
\\
\it (v) \rm\quad\quad &\text{There exists  }\kappa_0>0\text{
such that for }\rho_0:=1-\kappa_0a^{\frac{1}{n+1}}\text{ and $\beta_{\rho_0} := l_R(B_{\rho_0})$ we have }\\
&\quad\quad\text{(a)}\qd\qd H^{n-1}\lt(w\lt(\partial U_1\rt)\cap
\beta_{\rho_0} \rt)\leq c a,\nn\\
&\qd\qd\text{(b)}\qd\qd\deg( w, B_1, \xi) = 1\text{ for all }\xi\in \beta_{\rho_0} .
\end{align*}
\end{a5}
The ellipsoid  $\beta_{\rho} := l_R(B_\rho)$ should be thought of as the counterpart in the image $u(B_1)$
of the ball $B_\rho$ in the reference configuration.

\begin{proof}
We will apply a general truncation result, Lemma \ref{L7}, which is
proved in an Appendix. Toward this end, will write $f(Du) =
d^{s}(Du, K)$, where $s = 1+ \frac{p}{q^*}$ and $q^*$ is the Holder
conjugate of $q$. We first claim that
\begin{equation}
\label{L9.c1}
\| f(Du) \|_{W^{1,1}(B_1)} \leq  c a.
\end{equation}
We have already proved in (\ref{youngsineq}) that $\| D( f(Du))\|_{L^1(B_1)}\leq a$, so we only need to estimate
$\| f(Du)\|_{L^1}$. To do this, let us temporarily write $g(x) = d(Du\lt(x\rt),K)$, so that $f(Du) = g^s$.
Then (assuming $\ep$ is small enough) \eqref{L9.c1} follows from
\eqref{e1}, \eqref{youngsineq}, and the inequality
\begin{equation}
\| g^s \|_{L^1(B_1)} \leq C\left( \| g\|_{L^1(B_1)}^s + \| D (g^s)\|_{L^1(B_1)}\right)
\label{interp.sob}\end{equation}
since the terms on the right-hand side are just $\| d(Du,K)\|_{L^1}^s$ and $\| D(f(Du))\|_{L^1}$.

To prove \eqref{interp.sob},
we use Holder's inequality and the Sobolev embedding theorem to deduce that
\begin{align*}
\| g^s \|_{L^1\lt(B_1\rt)}  \ = \ \| g \|_{L^s\lt(B_1\rt)}^s
&\leq C \| g \|_{L^1\lt(B_1\rt)}^{s\theta} \ \| g \|_{L^{ns/(n-1)}\lt(B_1\rt)}^{s(1-\theta)}
\quad\quad\quad\quad\quad
\mbox{ for some $\theta\in (0,1)$}
\\
&= C   \| g \|_{L^1\lt(B_1\rt)}^{s\theta} \ \| g^s \|_{L^{n/n-1}\lt(B_1\rt)}^{(1-\theta)}\\
&\leq  C  \| g \|_{L^1\lt(B_1\rt)}^{s\theta}\lt(\| g^s \|_{L^1\lt(B_1\rt)}
+\| D( g^s)\|_{L^1\lt(B_1\rt)}\rt)^{(1-\theta)}.
\end{align*}
Then Young's inequality $ab \leq  \theta a^{1/\theta} + (1-\theta) b^{1/(1-\theta)}$ implies that
\[
\| g^s \|_{L^1\lt(B_1\rt)} \leq C \theta \|g \|_{L^1\lt(B_1\rt)}^s + (1-\theta)( \| g^s \|_{L^1\lt(B_1\rt)} + \| D(g^s)\|_{L^1\lt(B_1\rt)})
\]
which proves \eqref{interp.sob}.

We now fix $\lambda = \lambda(K)$ large enough that $d(Du, K) \geq
\frac 12 |Du|$ whenever $|Du|\geq \lambda$, and we apply Lemma
\ref{L7} to $u$ with this choice of $\lambda$, and with $q=1$ and
$f$ as above, so that $f(Du) \geq |Du| - c$. This produces a
Lipschitz function $w:B_1\to \R^n$ with $\|Dw\|_{L^\infty(B_1)}
\leq C \lambda  = C$. From conclusion {\em (ii)} of Lemma \ref{L7},
\begin{equation}
\label{xy187} \|Du-Dw\|_{L^1(B_1)} \leq \frac C \lambda \int_{\{
x\in B_1: |Du(x)|>\lambda \}} |Du| \overset{(\ref{e1})}{\leq} C \apep.
\end{equation}
Since $d(Dw, K) \leq d(Du,K) + |Du-Dw|$, it follows that
$\int_{B_1} d(Dw,K) \le C \ep$. Thus we have proved {\em (i)}.

Now let $U_0$ be the set constructed in
Lemma \ref{L3}. Recall that $d(Du,SO(n)A_i) = d(Du,K)$ in
$U_0$, so that
\begin{eqnarray}
\int_{B_1} d\lt(Dw,SO\lt(n\rt) A_i \rt)&\leq&\int_{U_0}
d\lt(Du,K\rt)+\lt|Du-Dw\rt| dx
+cL^n\lt(B_1\backslash U_0\rt)\nn\\
&\overset{(\ref{e1}),(\ref{L3.c1}),\eqref{xy187}}{\leq}&C a^{\frac{n}{n-1}}.\nn
\end{eqnarray}
So by the one-well $L^1$ Liouville Theorem\begin{footnote}{Strictly speaking,
we should change variables in a way that transforms $SO(n)A_i$ into $SO(n)$,
apply the $L^1$ Theorem, and then change variables  back; this is justified, since
the theorem we are citing is valid on any Lipschitz domain.}\end{footnote} (\cite{contsw} Proposition 2.6) there exists $R\in
SO\lt(n\rt)A_i$ such that $\int_{B_1} \lt|Dw-R\rt|\leq Ca$.  And by
Poincar\'e there exists an affine function $l_R$ with $Dl_R=R$ such
that $\int_{B_1} \lt|w-l_R\rt|\leq Ca$. By an interpolation
inequality, Theorem 5.9 \cite{adams}, this gives
$$
\|w-l_R\|_{L^{\infty}\lt(B_1\rt)} \leq C
\|w-l_R\|^{\frac{1}{n+1}}_{L^1\lt(B_1\rt)}
\|w-l_R\|^{\frac{n}{n+1}}_{W^{1,\infty}\lt(B_1\rt)} \leq C
a^{\frac{1}{n+1}}(\|w-l_R\|_{L^{\infty}\lt(B_1\rt)} +C)^{\frac{n}{n+1}},
$$
and this is easily seen to imply {\em (iv)}.

Next, Lemma \ref{L7} also asserts that there exists an open set
$E'\subset B_1$ with smooth boundary, such that $E :=  \{ x\in B_1:
u(x)\neq w(x)\} \subset E'$, and
\[
L^n(E')^{\frac{n-1}n} +  Per_{B_1}(E') \leq \|f(Du)\|_{W^{1,1}}  \overset{\eqref{L9.c1}}{\leq} Ca.
\]
We define $U_1 := U_0\setminus E'$. Then conclusion \it (ii) \rm is
immediate and conclusion \it (iii) \rm follows directly from the above
estimates of $E'$ and corresponding properties of $U_0$ from Lemma
\ref{L3}.

We now fix a constant $\kappa_0>0$ such that, if we define $\rho_0 := 1 - \kappa_0 a^{\frac 1{n+1}}$,
then
\begin{equation}
\beta_{\rho_0} = l_R(B_{\rho_0})  \subset w(B_1)\setminus  w(\partial B_1) .
\label{c1.def}\end{equation}
This is possible due to conclusion {\em (iv)}.
It follows that
\[
H^{n-1}(w(\partial U_1)  \cap \beta_{\rho_0} ) \leq H^{n-1}\left( w\lt(\partial U_1 \cap  B_1\rt)\right) \leq Ca
\]
using conclusion \it (iii), \rm together with the fact that $w$ is
Lipschitz, so we have shown \it (v), (a). \rm Finally, for $t\in [0,1]$ and $x\in B_1$, define $w_t(x) = t w(x) + (1-t)l_R(x)$.
It follows from conclusion {\em (iv)} that, taking $\kappa_0$ larger if necessary,
$w_t(\partial B_1)\cap \beta_{\rho_0} = \emptyset$ for every $t\in [0,1]$. Thus the
homotopy invariance of degree implies that for $\xi \in \beta_{\rho_0}$,
\[
\deg(w, B_1, \xi) =
\deg(w_1, B_1, \xi) =
\deg(w_0, B_1, \xi) =
\deg(l_R, B_1, \xi) = 1.
\]
This completes the proof of {\em (v)}.
\end{proof}


%
%
%
%

\section{Partial invertibility of $u$}\label{S4}

The main result of this section is the following
\begin{a5}\label{P.invert}
Suppose $u:B_1\subset \R^n\to \R^n$ is a smooth function satisfying
\eqref{e1}, with the wells labeled as in \eqref{e3} so that $A_1$ is
the majority phase, and let $U_1\subset B_1$ be the set found in
Proposition \ref{L9}. Then there exists a smooth  open set
$\UI_1\subset \R^n$ and a $C^1$ function $v:\UI_1\to U_1$ such that
$\| Dv\|_{L^\infty(\UI_1)}\le C$,
\begin{equation}
\label{we43}
u(v(\xi)) = \xi\text{ for all }\xi\in \UI_1,
\end{equation}
and for $\ep=\vs^{\frac{1}{p}}$
\begin{equation}
\int_{\UI_1} d(Dv, A_1^{-1}SO(n))\ dx \le \ C  \apep
\label{v.bound}\end{equation}
Moreover, there is an affine map $l_R$ with $Dl_R = R\in SO(n)A_1$  and a constant $\kappa_1$ such that
for $\rho_1 := 1 - \kappa_1 a^{\frac 1{n+1}}$and $\beta_{\rho_1}:= l_R(B_{\rho_1})$, the following hold:
\begin{equation}
\UI_1\subset \beta_{\rho_1}
\quad
Per_{ \beta_{\rho_1}}\UI_1 \le Ca,
\quad
L^n( \beta_{\rho_1}
\setminus \UI_1) \le C a^{\frac n{n-1}}
\label{UI.properties}\end{equation}
Finally, there exists a constant $k_1> \kappa_1$ such that
for $r_1 = 1-k_1 a^{\frac 1{n+1}}$,
\begin{equation}
L^n( B_{r_1} \setminus v(\UI_1)) \le C a^{\frac n{n-1}}.
\label{image.v}\end{equation}
\end{a5}

We start by proving a lemma in which
we find a set $D\subset B_1$ of small measure, and a radius $r_0$ close to $1$,
such that the Lipschitz approximation $w$ found earlier is one-to-one on $B_{r_0}\setminus D$.
We do not however have any information about the perimeter of $D$.
The properties \eqref{xy65a}, \eqref{xy65} of $D$ that we record in the statement of the Lemma
are consequences of these facts that will be useful in our later arguments.
The proof is follows arguments from \cite{contsw}.

\begin{a1}
\label{L12.5} Let $\beta_{\rho_0}$ be as defined in Proposition \ref{L9} (v).
Suppose $w:B_1\to \R^n$ is a Lipschitz function that satisfies (i), (iv) and (v) (b) of
Proposition \ref{L9}. Let  $V:=\lt\{x\in B_1:det\lt(Dw\lt(x\rt)\rt)\leq \sigma \rt\}$ and define
\begin{equation}
\label{xy43}
D:=\lt\{x\in B_1:w\lt(x\rt)\in w\lt(V\rt)\cap \beta_{\rho_0} \rt\}.
\end{equation}
Then $L^n\lt(D\rt)< C \apep$. In addition
there exists a constant $k_0>0$ such that for
$r_0:=1-k_0 a^{\frac{1}{n+1}}$, function $w$ is one-to-one on
$B_{r_0}\backslash D$ and for any set $S\subset B_{r_0}$,
\begin{equation}
\label{xy65a}
\deg(w, S, \xi) = 1
\quad\mbox{ for }\xi\in w(S\setminus D) = w(S)\setminus w(V);
\end{equation}
and
\begin{equation}
\label{xy65}
S\backslash D = w^{-1}\lt(w\lt(S\rt)\rt)\backslash D.
\end{equation}
\end{a1}

This lemma will be used not only in the proof of Proposition \ref{P.invert},
but also  in the proof of Lemma \ref{LL2}.

\begin{proof}
Recall in Proposition \ref{L9} we defined $\rho_0:=1-\kappa_0
a^{\frac{1}{n+1}}$, where the constant $\kappa_0$ was chosen so that
$\deg(w, B_1, \xi) =  1$ for all $\xi\in \beta_{\rho_0}= l_R(B_{\rho_0})$.

We choose $k_0$ so large that $w(B_{r_0}) \subset \beta_{\rho_0}$.
This is possible by Proposition \ref{L9}, {\em (iv)}.

{\em Step 1}. We first verify \eqref{xy65a} and \eqref{xy65}. If $S\subset B_{r_0}$, then
the choice of $r_0$ and the fact that $w$ satisfies {\em (v) (b)}
of Proposition \ref{L9}  imply that
$\mathrm{deg}\lt(w,B_1,\xi\rt)=1$ for any $\xi\in w\lt(S\rt)$. Thus for
$\xi\in w\lt(S\rt)\backslash w\lt(V\rt)$ we have
$$
\ca{w^{-1}\lt(\xi\rt)}=\sum_{x\in w^{-1}\lt(\xi\rt)} \mathrm{sgn}\lt(\det\lt(Dw\lt(x \rt)\rt)\rt)=1,
$$
Now  if $x\in w^{-1}\lt(w\lt(S\rt)\rt)\backslash D$, then
$w\lt(x\rt)\in w\lt(S\rt)\backslash w\lt(V\rt)$. Hence $w(x)$ has a unique preimage, which
necessarily belongs to $S$.
Thus
\[
w^{-1}\lt(w\lt(S\rt)\rt)\backslash D \subset S\backslash D.
\]
The opposite inclusion is obvious, so we have proved  \eqref{xy65}.
Similarly, if $\xi\in w(S\setminus D)$ then $w^{-1}(\xi)$ consists of
one point,
say $x$, which evidently belongs to $S\setminus D \subset S\setminus V$.
This implies that $\det Dw(x)>\sigma$. Thus
$\mathrm{sgn}\lt(\det\lt(Dw\lt(x\rt)\rt)\rt) = 1$, and \eqref{xy65a} follows.

\em Step 2. \rm We next show that $L^n\lt(D\rt)\leq C \apep$.

Note that from the choice \eqref{sigma.c2} of $\sigma$
and since $w$ satisfies \it (i) \rm of Proposition \ref{L9} we know that $L^n\lt(V\rt)\leq C \apep$.

Recall the change of variables degree formula (see \cite{fon} for
example)
\begin{equation}
\int_{\R^n} \psi\lt(\xi\rt)\, deg\lt(w,A,\xi\rt) d\xi
=\int_{A} \psi\lt(w\lt(x\rt)\rt)\det\lt(Dw\lt(x\rt)\rt) dx
\label{c.v.degree}\end{equation}
for open $A\subset B_1$ and $\psi\in L^\infty( \R^n)$.
We define
$\psi\lt(\xi\rt)=\chara_{w\lt(V\rt)\cap
\beta_{\rho_0}}\lt(\xi\rt)$, so that the definition \eqref{xy43} of $D$ implies
that $\psi(w(x)) = \chara_D(x)$. Then \eqref{c.v.degree}
yields
\begin{eqnarray}
\label{xy427.7}
\int_{w\lt(V\rt)\cap
\beta_{\rho_0}}deg\lt(w,B_1\backslash V,\xi\rt)
d\xi&=&\int_{D\cap(B_1\backslash V)}
 \det\lt(Dw\lt(x\rt)\rt) dx\nn\\
&\geq& \sigma\  L^n\lt(D\backslash V\rt).
\end{eqnarray}
Recall from  Proposition \ref{L9} (v) that
$deg\lt(w,B_1,\xi\rt)=1$ for all $\xi\in \beta_{\rho_0}$. Thus
\begin{eqnarray}
deg\lt(w,B_1\backslash V,\xi\rt)&\leq&1+\lt|deg\lt(w,V,\xi\rt)\rt|\nn\\
&\leq&1+\ca{w^{-1}\lt(\xi\rt)\cap V}.\nn
\end{eqnarray}
Note that
$\int_{w\lt(V\rt)} \ca{w^{-1}\lt(\xi\rt)\cap V} d\xi=\int_{V} \lt|\det\lt(Dw\lt(x\rt)\rt)\rt| dx\leq
\| Dw\|_{L^\infty}^n L^n(V) \le C \apep$. Similarly,
$\int_{w(V)} d\xi = L^n(w(V)) \le \|Dw\|_\infty^n L^n(V) \le C \apep$. Thus
\begin{eqnarray}
\int_{w\lt(V\rt)\cap \beta_{\rho_0}} deg\lt(w,B_1\backslash V,\xi\rt) d\xi &\leq&
\int_{w\lt(V\rt)\cap \beta_{\rho_0}} 1+\ca{w^{-1}\lt(\xi\rt)\cap V}
d\xi\nn\\
&\leq&C \apep.\nn
\end{eqnarray}
Putting this together with (\ref{xy427.7}) we have
$L^n\lt(D\backslash V\rt)\leq C \apep$. Since we know
$L^n\lt(D\cap V\rt)\leq C \apep$ this establishes $L^n\lt(D\rt)\leq C \apep$.
\end{proof}

\begin{a1}
\label{L12.5a} Let
$w:B_1\to \R^n$ be a Lipschitz function
satisfying the conclusions of Proposition \ref{L9}.
If $S\subset U_1\cap B_{r_0}$, for $r_0$ as defined in Lemma \ref{L12.5}, then
\begin{equation}
\label{eq21}
\lt|L^n\lt(w\lt(S\rt)\rt) - \det A_1 L^n\lt(S \rt)\rt|\leq  \apep.
\end{equation}
\end{a1}

\begin{proof}
For
$S\subset U_1\cap B_{r_0}$
it follows from Proposition \ref{L9} {\em (i), (ii)}
that $\int_S d(Dw, SO(n)A_1) \le C \apep$.
Then using  
the fact (Proposition \ref{L9} {\em (v)}) that
$\deg(w, B_1,\xi)=1$ for every $\xi\in w(S)\subset \beta_{\rho_0}$, we
compute
\begin{align*}
L^n(w(S))
&=
\int_{\R^n} \chara_{ w(S)}(\xi) \  \deg(w, B_1,\xi)\ d\xi  \\
&\overset{\eqref{c.v.degree}}=
\int_{ w^{-1}(w(S))} \det Dw(x)\ dx\\
&=
\int_{ w^{-1}(w(S)) \setminus D} \det Dw(x)\ dx
+
\int_{ w^{-1}(w(S)) \cap D} \det Dw(x)\ dx
\end{align*}
where $D$ was defined in the previous lemma, in which we also proved that
$L^n(D)\le C\apep$.
To estimate the second term, note that $|\int_{ w^{-1}(w(S)) \cap D} \det Dw(x)\ dx| \le C L^n(D) \le C\apep$.
And in view of \eqref{xy65},
\begin{align*}
&\int_{ w^{-1}(w(S)) \setminus D} \det Dw(x)\ dx
=
\int_{ S \setminus D}  \det Dw(x)\ dx\\
&\hspace{5em}
= \det A_1 [L^n(S) -
L^n(S\cap D) ]
+
\int_{ S \setminus D} ( \det Dw(x)- \det A_1)\ dx.
\end{align*}
Since $Dw$ is Lipschitz, $\int_S |\det Dw-\det A_1| \le C \int_S d(Dw, SO(n)A_1) \le  C \apep$. So combining the above inequalities and
using again the fact that $L^n(D)\le C \apep$, we obtain
\eqref{eq21}.
\end{proof}

%
%
%
%
%

\begin{proof}[Proof of Proposition \ref{P.invert}]
Throughout the proof we will use notation introduced in Proposition \ref{L9}.
Recall also that  $r_0 = 1 - k_0 a^{\frac 1{n+1}}$ was fixed in Lemma \ref{L12.5}.
We fix $\rho_1= 1 - \kappa_1 a^{\frac 1{n+1}}$ by choosing
$\kappa_1>k_0$ large enough that
$l_R(B_{\rho_1}) = \beta_{\rho_1}\subset w( B_{r_0})\setminus
w(\partial B_{r_0})$. This is possible
due to Proposition \ref{L9} {\em (iv)}.

Next, we define
\begin{equation}
\label{weib2}
\UI_1 := \{ \xi\in \beta_{\rho_1}\ : \ \deg(w, U_1\cap B_{r_0}, \xi) = 1\}.
\end{equation}
{\em Step 1}. First we establish some properties of $\UI_1$.
General facts about degree imply that $ \deg(w, U_1\cap B_{r_0},  \cdot)$ is locally constant
in $\R^n\setminus w( \partial(U_1\cap B_{r_0}))$. Since $w( \partial(U_1\cap B_{r_0}))$ is closed,
it follows that $\UI_1$ is open. In addition, we deduce that
\[
\partial \UI_1 \cap \beta_{\rho_1}\subset w( B_{r_0}\cap \partial U_1).
\]
Since $w$ is Lipschitz, it follows from conclusion {\em (v)} of Proposition \ref{L9}
that
\begin{equation}
{Per}_{\beta_{\rho_1}}\lt(\UI_1\rt)
\le
H^{n-1}(w( B_{r_0}\cap \partial U_1))\ \le C H^{n-1}( B_{r_0}\cap \partial U_1)
\le Ca.
\label{UI.perest}\end{equation}
Next,  recall that Lemma \ref{L12.5} implies that
$\deg(w, U_1\cap B_{r_0}, \xi) = 1$ for every $\xi\in w( ( U_1\cap B_{r_0})\setminus D)$,
see \eqref{xy65a}, where $D\subset B_1$ is defined in \eqref{xy43} and has the property that
$L^n(D)\le C \apep$.
Again using Proposition \ref{L9} {\em (iv)}, we know that $w(B_{1/2}) \subset \beta_{\rho_1}$
if $a$ is small enough. For such $a$,
it follows that $w(( U_1\cap B_{1/2})\setminus D)\subset \UI_1$,
and Lemma \ref{L12.5a} with Proposition \ref{L9} {\em (i)}, {\em  (ii)} implies that
\[
L^n(\UI_1) \ge L^n(w(( U_1\cap B_{\frac{1}{2}})\setminus D)\,) \ge \ \det A_1 \
L^n( ( U_1\cap B_{\frac{1}{2}})\setminus D) - C \apep.
\]
Then Proposition \ref{L9} {\em (iii)} and the fact that $L^n(D)\le C \apep$
yield
$L^n(\UI_1) \ge  \det A_1 L^n( B_{1/2} ) - C a^{\frac n{n-1}}- C \apep$.

On the other hand, we know from \eqref{UI.perest} and the relative isoperimetric inequality that
\[
\min \{ L^n(\UI_1),  L^n (\beta_{\rho_1}\setminus \UI_1) \} \le C a^{\frac n{n-1}}
\]
Combining these facts, we conclude that
$ L^n (\beta_{\rho_1}\setminus \UI_1) \le C a^{\frac n{n-1}}$.

Thus we have verified that $\UI_1$ has all the properties asserted in \eqref{UI.properties}.

{\em Step 2}.
For $\xi\in \UI_1$,
recalling that $\det Dw>0$ in $U_1$,
we deduce from the definition of $\UI_1$  that
\begin{align*}
1 = \deg(w, U_1\cap B_{r_0}, \xi)
&=
\sum_{\{ x\in U_1\cap B_{r_0}:  w(x) = \xi \} }\mathrm{sgn}\lt(\det\lt(Dw\lt(x\rt)\rt)\rt) \\
&= \mbox{Card}(\{ x\in U_1\cap B_{r_0}:  w(x) = \xi \} ).
\end{align*}
It follows that
it makes sense to define $v:\UI_1\to \R^n$ by stipulating that
\begin{equation}
v(\xi) = x \iff
x\in U_1\cap B_{r_0}\mbox{ and }w(x)=\xi.
\label{v.def}\end{equation}
Since $u=w$ in $U_1$, we deduce that $u(v(\xi)) = \xi$ for all $\xi\in \UI_1$, as required.

We next verify that $v$ is $C^1$.
To do this, fix any $\xi\in \UI_1$, and let $x=v(\xi)$.
Since $w$ is smooth with $\det Dw \ne 0$ in $U_1$, the inverse function theorem
implies that there is a neighborhood $N_\xi$ of $\xi$ and a $C^1$
map $\tilde v: N_\xi\to B_1$ such that $\tilde v(\xi) = x$ and
$w(\tilde v(\eta)) = \eta$ for all $\eta$ in $N_\xi$. Since $U_1$ is open,
we may assume (by shrinking $N_\xi$ if necessary) that
$\tilde v(\eta)\in U_1\cap B_{r_0}$ for all $\eta\in N_\xi$.
Then it is clear that $\tilde v = v$ in $N_\xi$, and therefore that
$v$ is $C^1$.

{\em Step 3}.
Since $v$ is $C^1$, it follows that $Dv(\xi) = Dw(v(\xi))^{-1}$.
Thus by the change of variables $\xi = w(x)$
(which is straightforward, since $v$ is
a bijection onto its image) we find that
\[
\int_{\UI_1} d(Dv(\xi), A_1^{-1} SO(n)) \ d\xi =
\int_{v(\UI_1)} d(Dw(x)^{-1}, A^{-1}SO(n)) \det Dw(x) \ dx.
\]
Because from (\ref{v.def}) we know $v\lt(\UI_1\rt)\subset U_1\cap B_{r_0}$ and
$d(Dw(x), SO(n)A_1) \le  \sigma$ for $x\in U_1\subset U_0$, so
there is a constant $C$ such that $\det Dw(x) \le C$ and
\[
d(Dw(x)^{-1}, A_1^{-1}SO(n)) \le C \ d( Dw(x), SO(n)A_1) = C \ d(Dw(x), K)
\]
for all $x\in U_1$. We conclude that
\[
\int_{\UI_1} d(Dv(\xi), A_1^{-1}SO(n)) \ d\xi \le C \int_{B_1} d(Du,K) \overset{(\ref{e1})}{\le} C \apep.
\]

{\em Step 4}. Finally, we fix $k_1$ such that, if we define $r_1 = 1- k_1 a^{\frac 1{n+1}}$,
then $w(B_{r_1}) \subset \beta_{\rho_1}$. This is possible as usual due to Proposition \ref{L9} {\em (iv)}.
The definitions imply that $r_1<r_0$.

Then the definition \eqref{v.def} of $v$ implies that
$B_{r_1}\setminus v(\UI_1) = (B_{r_1}\setminus U_1) \cup S$,
where
\[
S := \{ x\in B_{r_1}\cap U_1 \ : \ w(x)\not\in \UI_1\}.
\]
Then we can use Proposition \ref{L9} {\em (iv)}
and Lemma \ref{L12.5a} to estimate
\begin{align*}
L^n( B_{r_1}\setminus v(\UI_1))
&\le
L^n( B_{r_1}\setminus U_1 ) + L^n(S) \\
&\le
C a^{\frac n{n-1}} +\frac{L^n(w(S))}{\det A_1} + C \ep.
\end{align*}
And $w(S)\subset \beta_{\rho_1}\setminus \UI_1$, so \eqref{UI.properties} implies
that $L^n(w(S)) \le C a^{\frac n{n-1}}$. This proves \eqref{image.v}
and completes the proof of the Proposition.
\end{proof}

\section{Non-shrinking pairs}\label{S5}

In this section we prove

\begin{a5}
\label{L10} Suppose that $u:B_1\subset \R^n\to \R^n$ is a smooth function satisfying (\ref{e1})
and assume as in (\ref{e3}) that $SO(n)A_1$ is the majority phase. Then there exists $\GI\subset B_{r_1}\times B_{r_1}$
(where $r_1  = 1 - k_1 a^{\frac 1{n+1}}$ was fixed in Proposition \ref{P.invert} and satisfies \eqref{image.v})
such that for $\ep=\vs^{\frac{1}{p}}$,
\begin{equation}
\mbox{ if }(x,y)\in \GI, \mbox{ then } \
\large\left| |u\lt(y\rt)-u\lt(x\rt)| \  -  \  |A_1(x-y)| \large\right| \le  C  \ep
\label{L10c}\end{equation}
and
$L^{2n}( ( B_{r_1}\times B_{r_1})\setminus \GI) \le C a^{\frac 1p}$.
\end{a5}

Following the proof we give a couple of corollaries that will be used in later
sections. Throughout the proof we will use notation from the statement of
Proposition \ref{P.invert}.

\begin{proof}[ Proof of Proposition \ref{L10}. ]
It suffices to find a set $\GI_2\subset B_{r_1}\times  B_{r_1}$ such that
\begin{equation}
L^{2n}( ( B_{r_1}\times B_{r_1})\setminus \GI_2) \le C a^{\frac 1p}
\label{G2.1}\end{equation}
and
\begin{equation}
\mbox{ if }(x,y)\in \GI_2, \mbox{ then } \
 |u\lt(y\rt)-u\lt(x\rt)|  \geq  \  |A_1(x-y)| - C \ep
\label{G2.2}\end{equation}
since then the conclusions of the lemma follow if we define $\GI := \GI_1\cap \GI_2$,
where $\GI_1$ was constructed in Lemma \ref{L6}.

{\em Step 1} .
We define
\begin{equation}
\Gamma :=
\{ (\xi,\eta)\in  \beta_{\rho_1}\times  \beta_{\rho_1}\ : \ [\xi,\eta] \subset \UI_1, \int_{[\xi,\eta]} d(Dv, A_1^{-1}SO(n)) \ dH^1 \ \le \ep\}
\label{Gamma.def}\end{equation}
and
\begin{equation}
\GI_2 := \{ (v(\xi), v(\eta)) \ : \ (\xi,\eta) \in \Gamma\}.
\label{G2.def}\end{equation}
We claim that
\begin{equation}
L^{2n}( \left(\beta_{\rho_1}\times  \beta_{\rho_1}\right)\setminus \Gamma) \le Ca^{\frac 1p}.
\label{Gamma1}\end{equation}
and
\begin{equation}
\mbox{ if $(\xi,\eta)\in \Gamma$, then
$|A_1(v(\xi) - v(\eta))| \le |\xi-\eta|+ C \ep$}
\label{Gamma0}\end{equation}
In fact, \eqref{Gamma1} follows from {\em exactly}
the same argument used in the proof of Lemma \ref{L6} to establish \eqref{G1.1}.
That proof
relied on the facts that $L^n(B_1\setminus U_0) \le C a^{\frac n{n-1}}$, $H^{n-1}(\partial U_0\cap B_1) \le Ca$,
and $\int_{U_0} d(Du, SO(n)A_1)\le C \ep$.
In Proposition \ref{P.invert} we  proved that
the same estimates
hold with $B_1, U_0$ and $u$  replaced by $\beta_{\rho_1}, \UI_1$ and $v$,
respectively, so the earlier
arguments can be repeated  word for word.
Next fix $(\xi,\eta)\in \Gamma$, and write $\tau := \frac{\eta-\xi}{|\eta-\xi|}$. Then
\begin{align*}
|A_1(v(\eta) - v(\xi))|
&=
\left| \int_{[\xi,\eta]}  A_1 Dv \ \tau \  dH^1 \right| \\
&\le
\int_{[\xi,\eta]} \left[1 + C \,d(Dv,A_1^{-1}SO(n)) \right] dH^1
\ \overset{(\ref{Gamma.def})}{\le} \
|\eta-\xi| + C \ep
\end{align*}
for $(\xi,\eta)\in \GI_1$, proving \eqref{Gamma0}.

{\em Step 2}.
Observe that \eqref{G2.2} is an immediate consequence of \eqref{Gamma0} and the definition \eqref{G2.def}
of  $\GI_2$ and equation (\ref{we43}) of Proposition \ref{P.invert}. To verify \eqref{G2.1},
note that
\[
(B_{r_1} \times B_{r_1}) \setminus \GI_2 \ \  \subset \ \
\left[(B_{r_1} \times B_{r_1}) \setminus (v(\UI_1)\times v(\UI_1))\right]
\ \cup
\ \SI
\]
for $\SI := \{ (x,y)\in v(\UI_1)\times v(\UI_1) \ : \ (u(x), u(y))\not\in \Gamma\}$.
We deduce from \eqref{image.v} that
\[
L^{2n}( \ (B_{r_1} \times B_{r_1}) \setminus (v(\UI_1)\times v(\UI_1)) \ ) \le C a^{\frac n{n-1}}.
\]
To estimate the measure of $\SI$, we use
a change of variables (which is straightforward, since $v$ is a $C^1$ diffeomorphism onto its image) to compute
\begin{align*}
L^{2n}(\SI)
&=  \
\int_{v(\UI_1)\times v(\UI_1) }
\chara_{ (u(x), u(y))\not\in \Gamma }  \ dx \ dy\\
&= \
\int_{\UI_1\times \UI_1} \
\chara_{(\xi,\eta)\not \in \Gamma} \  \det Dv(\xi) \ \det Dv(\eta) \ d\xi \ d\eta \\
&\le \
J^2 \ L^{2n} ( (\UI_1\times \UI_1)\setminus \Gamma) + \int_{\UI_1\times \UI_1}
|J^2 - \det Dv(\xi) \ \det Dv(\eta)| \ d\xi \ d\eta
\end{align*}
for $J:= \det A_1^{-1}$.
The integral on the right-hand side is bounded above by
\[
\int_{ \UI_1\times \UI_1 }
J |J - \det Dv(\xi)| + |\det Dv(\xi)||J  -  \ \det Dv(\eta)| \ d\xi \ d\eta
\overset{\eqref{v.bound}}{ \le} C  \apep
\]
since $|J - \det Dv(\xi)| \le C d(Dv(\xi), A^{-1} SO(n))$ in $\UI_1$  (because
$v$ is Lipschitz). And Step 1 (recall definition (\ref{weib2})) implies that
\[
L^{2n} ( (\UI_1\times \UI_1)\setminus \Gamma)  \ \overset{\eqref{weib2}}{\le} \
L^{2n} ( (\beta_{\rho_1}\times \beta_{\rho_1}) \setminus \Gamma)   \overset{\eqref{Gamma1}}{\le} Ca^{1/p}.
\]
Combining the above estimates, we arrive at \eqref{G2.1}.
\end{proof}

Our first Corollary is

\begin{a6}
Assume the hypotheses of Proposition \ref{L10}. Then
for any $C_1>0$, the set of points $x\in B_{r_1}$ such that
\begin{equation}
L^n(\{ y\in B_{r_1} \ : \
(x,y)\not\in \GI \})  \ge C_1 a^{1/p}
\label{L10.cor1.c}\end{equation}
has measure at most $\frac C {C_1}$, where $\GI$ is as found  in Proposition \ref{L10}.
\label{L10.cor1}
\end{a6}

\begin{proof}
This follows from Fubini's Theorem, Chebyshev's inequality, and the
conclusion of Proposition \ref{L10}, ie the  fact that $L^{2n}(B_{r_1}\times B_{r_1}\setminus \GI)\le Ca^{1/p}$.
\end{proof}

\begin{a6}
Assume the hypotheses of Proposition \ref{L10}.

Suppose also that $B_\delta(y_0),\ldots, B_\delta(y_n)$ are pairwise disjoint balls contained
in $B_{r_1}$, and that $\HI\subset B_1$ is a measurable set such that $L^n(B_1\setminus\HI)\leq
\frac{1}{2} L^n(B_\delta)$.

Then for every sufficiently large $C_1>0$, there exists $a_0>0$
such that if $a\le a_0$, then there are
points $x_k\in B_\delta(y_k)\cap \HI$ for $k=0,\ldots, n$ such that
\begin{equation}
\large\left| \ |u\lt(x_k\rt)-u\lt(x_l\rt)| \  -  \  |A_1(x_k - x_l)|  \ \large\right| \le  C \ep
\quad\mbox{ for all }k\ne l
\label{rsmp1}\end{equation}
(for the same $C$ as in \eqref{L10c}),  and
\begin{equation}
L^n(\{ y\in B_{r_1} \ : \
(x_k,y)\not \in \GI \})  \le C_1 a^{1/p} \quad\mbox{ for all }k
\label{rsmp2}\end{equation}
\label{C.rsimp}
\end{a6}

We prove in Lemma \ref{L12}, in an appendix, that if \eqref{rsmp1} holds, and if $\{x_0,\ldots, x_n\}$
are the vertices of a nondegenerate simplex, then there exists an affine function $l_R$
with $Dl_R = R\in SO(n)A_1$ such that $| u(x_i) - l_R(x_i)| \le C \ep$ for every $i$.

\begin{proof}
Let us say that a point $x\in B_{r_1}$ is {\em good} if it does {\em not} satisfy \eqref{L10.cor1.c}, for some value $C_1>0$ to be selected below. Thus, for a fixed
good point $x$,  all $y\in B_{r_1}$ away from an exceptional set (depending on $x$)
of measure at most $C_1 a^{1/p}$ satisfy
$| \ |u(x)-u(y)| - |A_1(x-y)| \ | \le C\ep$.

Let us set $\HI' := \{ x\in \HI \ : \ x\mbox{ is good} \}$. We fix $C_1$ in the definition of ``good''
so large that
$L^n(B_{r_1}\setminus \HI') \le \frac{3}{4} L^n(B_\delta)$;
it follows from the hypothesis on $\HI$ and
from Corollary \ref{L10.cor1}
that this is possible.
We will show that if $a$ is small enough, then
there exist $x_k$, $k= 0,\ldots, n$, such that
\begin{equation}
\mbox{$x_k\in B_\delta(y_k)\cap \HI'$ and  \
$(x_k,x_l)\in \GI$ for every $k\ne l$.}
\label{xk.ind}\end{equation}
In view of the definitions, this will prove the corollary.

We fix $x_0$ be any point in $B_{\delta}(y_0)\cap \HI'$.
Now assume by induction that we have found $x_0,\ldots, x_{k-1}$ satisfying \eqref{xk.ind},
for some $k\le n$.
Since $x_0,\ldots, x_{k-1}$ belong to $\HI'$ and hence are good points, it follows that
\[
L^n\left( \bigcup_{l=0}^{k-1}\{ x\in B_\delta(y_k)  \cap \HI' \ : \  \mbox{$(x_l,x)\not\in \GI$}
\} \right) \le k C_1 a^{1/p}.
\]
Thus
\begin{align*}
&L^n( \{ x\in B_\delta(y_k)  \cap \HI' \ : \  \mbox{$(x_l,x)\in \GI$ for $l=0,\ldots, k-1 $}
\} ) \\
&\quad\quad \ge L^n(B_\delta(y_k)  \cap \HI' ) - k C_1 a^{1/p}
\ \ \ge  \ \frac 14 L^n(B_\delta) -  n C_1 a^{1/p}.
\end{align*}
In particular,  the above set is nonempty for every $k\le n$
if  $a \le a_0=\left[\frac 1{8 n C_1} L^n(B_\delta)\right]^p$.
Then we can pick $x_k$ to be any point in the above nonempty set, and
this eventually yields a collection satisfying \eqref{xk.ind}.
 \end{proof}

\section{Proof of Theorem \ref{T1}}\label{S6}

In this section we present the proof of Theorem  \ref{T1}.
Most of the work of the proof is carried out in two lemmas, in which we consider  the
case $\Omega = B_1$ and $\Omega' = B_\delta$ for some small $\delta>0$.
We start by assuming  these two lemmas
hold, and we use them to complete the proof of the theorem. The proofs of the lemmas follow.

\begin{proof}[Proof of Theorem \ref{T1}]

{\em Step 1}. If $m=2$, we claim that each well must satisfy either \eqref{T2.h1} or \eqref{T2.h2}.
To see this, suppose that $SO(n)A_1$ and $SO(n)A_2$ are two distinct wells such that
\eqref{T2.h1} does not hold for $i=1$.
We may assume by a polar decomposition that $A_1,A_2$ are both symmetric.
The assumption that \eqref{T2.h1} fails for $i=1$ says that
$|A_2  v|^2  \ge |A_1 v|^2$ for all $v$.
In particular, since this holds for $v$ of the form $v = A_2^{-1}w$, it follows that
\begin{equation}
\mbox{$|w|^2\ge    |A_1 A_2^{-1} w|^2 = w^T (A_2^{-1} A_1^2A_2^{-1} )w $ \ \  for all $w$,}
\label{1le2}\end{equation}
or in other words that $A_2^{-1} A_1^2A_2^{-1} \le Id$.
By taking inverses we find that $A_2 A_1^{-2}A_2 \ge Id$, or equivalently
$|w^T A_2 A_1^{-1}|^2 \ge |w|^2$ for all $w$. This in turn implies that
$| v^T A_1^{-1} |^2 \ge |v^T A_2^{-1}|^2$ for all $v$. To prove that \eqref{T2.h2} holds,
we must show that strict inequality holds for some $v$, which however is
clear, since otherwise equality would hold in \eqref{1le2}, which would imply
that $A_1 A_2^{-1}\in SO(n)$, and this is impossible since the two wells were assumed to be distinct.

Thus it suffices to show that the lemma holds if each well satisfies   \eqref{T2.h1} or \eqref{T2.h2}.

{\em Step 2}. Now fix a connected set $\Omega' \subset\subset \Omega$, and fix $r< \mbox{dist}(\Omega', \partial \Omega)$.
For $\delta$ to be specified below, we fix points $x_1,\ldots, x_N\in \Omega'$ such that $\Omega'\subset \cup_{k=1}^n B_{\delta r}(x_k)$.
For each $k$,
\begin{equation}
\frac 1\vs
\int_{B_r(x_k)} d^p\lt(Du,K\rt) + \frac { \lt|D^2 u\rt|^q}{\vs^q} dx
\leq
\frac 1\vs
\int_{\Omega} d^p\lt(Du,K\rt) + \frac { \lt|D^2 u\rt|^q}{\vs^q} dx
\leq a
\label{Brxk.est}\end{equation}
and so we can apply a suitable scaled version of Lemma \ref{L3} to find
some $i = i(k)$ and a set $U_0^k\subset B_r(x_k)$ satisfying
\eqref{L3.c1}, \eqref{a42} hold
(with $B_1$ and $A_i$ replaced by $B_r(x_k)$ and $A_{i(k)}$, and with the constants now
depending on $r$, which however has been fixed.)  These conclusions
imply that $i(k) = i(k')$ for any $k,k'$ such that $L^n( B_r(x_k) \cap B_r(x_k')) \ge C a^{\frac n{n-1}}$
for a suitable constant $C$. Thus by taking $a$ small enough we deduce that $i(k)$
is in fact independent of $k$, so that every ball $B_r(x_k)$ has the same majority phase.
We relabel the wells as usual so that $A_1$ represents this majority phase.

{\em Step 3}. By assumption $A_1$ satisfies \eqref{T2.h1} or \eqref{T2.h2}. In the
former case, it follows by continuity that $A_1$ satisfies the hypothesis
\eqref{hyp.rewrite} of Lemma \ref{LL1} (proved below) for some $\alpha>0$,
and similarly, if \eqref{T2.h2} holds, then hypothesis \eqref{hyp.rewrite2} of Lemma \ref{LL2} is valid for some $\alpha>0$.
We now require that $\delta\le \delta_i$, $i=1,2$, where $\delta_1, \delta_2$
appear in the conclusions of Lemmas \ref{LL1} and \ref{LL2} respectively.
Then  in view of \eqref{Brxk.est}, if $a$ is small enough then
we can apply Lemma \ref{LL1} or \ref{LL2} (scaled to a ball of radius $r$)
on each $B_r(x_k)$ to conclude that
\[
\int_{\Omega'} d\lt(Du,SO\lt(n\rt)A_1\rt)
\le C \sum_{k=1}^N\int_{B_{\delta r}(x_k)}
d\lt(Du,SO\lt(n\rt)A_1\rt)
\le C  \vs^{1/p}.
\]
Note also that  $d^p(Du, SO(n)A_1) \le  C [d(Du, SO(n)A_1) + d^p(Du,K)]$. Thus the above
inequalities immediately imply that $ \int_{\Omega'} d^p\lt(Du,SO\lt(n\rt)A_1\rt)  \le C  \vs^{1/p}$.
Finally, by applying
Theorem 3.1 of \cite{fmul} (see \eqref{mthm1}) we conclude that if $p>1$, then
\[
\inf_{R\in K} \| Du - R \|_{L^p(\Omega')}^p \le \vs^{1/p}.
\]
\end{proof}

%
%

We now give the proofs of the two lemmas used above. The first uses
Lemma  \ref{L.IbyP}, which is proved in the introduction.

\begin{a1}
\label{LL1} Let $\{A_1,\ldots, A_m\}$ be a set of $n\times n$
matrices, and let $K = \bigcup_i SO(n) A_i$. Let $u:B_1\to \R^n$ be
a smooth function that satisfies (\ref{e1}) and assume the matrices
have been labeled so that $A_1$ is the majority phase, i.e.\ the set
$U_0$ we obtain from Lemma \ref{L3} satisfies (\ref{e3}).

Suppose $A_1$ has the property that there exists $v_1\in S^{n-1}$ and $\alpha>0$ such
that
\begin{equation}
\mbox{ $\lt|A_1 v\rt|> (1+\alpha) \lt|A_j v\rt|\text{ for all }j\ge 2$ and all $v$ such that $|v\cdot v_1|> (1-\alpha)|v|$}.
\label{hyp.rewrite}
\end{equation}
Then there exist constants $a_0, \delta_1>0$ such that if $a\le a_0$ in \eqref{e1},
then
\begin{equation}
\label{wew55.56}
\int_{B_{\delta_1}} d\lt(Du,SO\lt(n\rt)A_1\rt) dx\leq C \ep.
\end{equation}
\end{a1}

\begin{proof}
{\em Step 1}.
We first find points $\{x_0,\ldots, x_n\}$
such that the hypotheses of Lemma \ref{L.IbyP} are satisfied,
together with some other conditions that we will need below.

Fix $0< \delta \le  \frac 18$, $\alpha>0$, and $y_0,\ldots, y_n$ such that
$|y_k| = \frac{1}{2}$ for all $k$, and
if $x_k\in B_\delta(y_k)$ for $k=0,\ldots, n$ and $x\in B_\delta(0)$, then
\begin{equation}
|\tau_k \cdot v_1| \ge 1- \alpha\quad\mbox{ for }\tau_k := \frac{x-x_k}{|x-x_k|}
\label{delta}\end{equation}
and
\begin{equation}
x = \sum_{k=0}^n\lambda_k x_k\mbox{ with }\sum_{k=0}^n\lambda_k=1 \ \ \mbox{ and } \ \
\lambda_k>\delta \mbox{ for all }k.
\label{delta2}
\end{equation}
For example, we can take $y_0 = \frac 12 v_1$, and $y_k = - s_1 v_1 + s_2 z_k$,
where $s_1, s_2$ are constants such that $s_1^2 + s_2^2 = \frac 14$, and
$\{z_1,\ldots, z_n\}$ are the vertices of a regular $n-1$-dimensional simplex sitting
on the unit  sphere in $v_1^\perp$. If $s_2$ is sufficiently small, then
$|\frac {y_k}{|y_k|} + v_1|<\alpha$, and  the above conditions hold if $\delta$
is sufficiently small.

We will write
\[
\HI := \{x\in B_1 \ : \ \int_{B_1} d(Du(z),K)|x-z|^{1-n} \ dz \ \le C_2 \ep\}.
\]
for a constant $C_2$ to be determined below.
Fubini's Theorem
implies that
\[
\int_{B_1} \int_{B_1} d(Du(z),K)|x-z|^{1-n} \ dz \ dx \ \le
C \int_{B_1} d(Du(z),K) \ dz \  \overset{\eqref{wieb1}}{\le}   C\ep
\]
so we deduce  from Chebyshev's inequality that $L^n(B_1\setminus \HI ) \le \frac C{C_2}$.
We now fix $C_2$ large enough that $L^n(B_1\setminus \HI) \le \frac 12L^n(B_\delta)$. Then Corollary
\ref{C.rsimp} implies that if $a$ is smaller than an appropriate constant $a_0$,
we can find points $x_k\in B_\delta(y_k)\cap \HI$ for $k=0,\ldots, n$
such that
\[
\large\left| \ |u\lt(x_k\rt)-u\lt(x_l\rt)| \  -  \  |A_1(x_k - x_l)|  \ \large\right|\le
C \ep
\quad\mbox{ for all }k\ne l.
\]
This is exactly the hypothesis of Lemma \ref{L12} (proved in Section \ref{SLA}).
The conclusion of this lemma is
that there exists an affine map $l_R$ with $Dl_R = R\in O(n)A_1$ such that
\begin{equation}
\mbox{$|u(x_k) - l_R(x_k)| <  C \ep$ for $k=0,\ldots, n$.}
\label{xy50a}\end{equation}

{\em Step 2}. It follows from \eqref{delta}, and \eqref{hyp.rewrite} that
for $\{x_0,\ldots\ x_n\}$ as found above, the hypotheses \eqref{shrink.dir} of Lemma \ref{L.IbyP} are satisfied
for every $x\in B_\delta$.
It follows from the lemma and  \eqref{xy50a} that
\begin{equation}
\label{wei5}
\sum_{k=0}^n \int_{[x_k, x]} d(Du, SO(n) A_1) \, dH^1 \ \le  \ C
\sum_{k=0}^n \int_{[x_k, x]} d(Du, K) \, dH^1 \ + C  \ep
\end{equation}
with a fixed constant $C$ valid for all $x\in B_\delta$.

{\em Step 3}.
To complete the proof we will integrate the above inequality over $x\in B_\delta$.
Both sides of the resulting estimate contain terms of the form
$\int_{B_\delta} F_k(x) \ dx$, for $F_k$ of the form $F_k(x) = \int_{[x_k,x]} f(y) dH^1$.
Note that by Lemma \ref{AX3}
\begin{align*}
\int_{B_\delta} F_k(x) \ dx
&=
\int_{\theta\in S^{n-1}} \int_0^\infty F_k(x_k + r\theta) r^{n-1} \chara_{ {x_k+r\theta\in B_\delta} } \ dr d H^{n-1}\theta\\
&=
\int_{\theta\in S^{n-1}} \int_0^\infty \int_0^r f(x_k + s\theta) r^{n-1} \chara_{ { x_0+r\theta\in B_\delta} } \ ds\ dr d H^{n-1}\theta.
\end{align*}
We apply Fubini's Theorem, integrate in the $r$ variable, and undo the transformation to polar coordinates,
to obtain
\begin{equation}
\label{wei2}
\int_{B_\delta} F_k(x) \ dx = \int f(x) G_k(x) |x-x_k|^{1-n} \ dx,
\quad
\mbox{ where }\ \
G_k(x_k + s\theta) := \int_s^\infty r^{n-1} \chara_{ { x_k+r\theta\in B_\delta} } dr.
\end{equation}
Now we integrate both sides of the inequality \eqref{wei5}. Since $G_k$ is clearly bounded,
\[
\sum_k
\int_{B_{\delta}}
(\int_{[x_k, x]} d(Du, K) \, dH^1  ) dx
 \le \sum_k C\int_{B_1} d(Du, K) |x-x_k|^{1-n}  \ dx  \ \le \  C \ep
\]
where we have used
the fact that $x_k\in \HI$ for $k=0,\ldots, n$. It is also easy to check that
$G_k(x)\,|x-x_k|^{1-n}\ge C^{-1}$ in $B_{\delta}$, which implies that
\[
\sum_k  \int_{B_{\delta}} d(Du, SO(n) A_1)
\le
C\sum_k
\int_{B_{\delta}}
(\int_{[x_k, x]} d(Du, SO(n) A_1) \, dH^1  ) dx.
\]
By combining these with \eqref{wei5} and defining $\delta_1 := \delta/2$, we complete the proof of the lemma.
\end{proof}

The proof of Theorem  \ref{T1} will be completed by the following lemma. As mentioned in
the introduction, the idea is roughly to apply to $u^{-1}$ an argument like that used in
the above lemma. Because $u$ is not invertible, we work with the Lipschitz approximation
$w$ found earlier, and we use a lemma (proved in Section  \ref{S:app}) that more or less allows
us to find a Lipschitz path in $B_1$ in the inverse image of a.e. line segment in
$\beta_{\rho_1}$.

\begin{a1}
\label{LL2} Let $\{A_1,\ldots, A_m\}$ be a set of $n\times n$
matrices, and let $K = \bigcup_i SO(n) A_i$. Let $u:B_1\to \R^n$ be
a smooth function that satisfies (\ref{e1}) and assume the matrices
have been labeled so that $A_1$ is the majority phase, i.e.\ the set
$U_0$ we obtain from Lemma \ref{L3} satisfies (\ref{e3}).

Suppose $A_1$ has the property that there exists $v\in S^{n-1}$ and $\alpha\in\lt(0,1\rt)$ such
that
\begin{equation}
\mbox{ $\lt|v^T A_1^{-1} \rt|> (1-\alpha)^{-1} \lt|v^T A_j^{-1} \rt|\text{ for all }j\ge 2$ and all $v$ such that $|v\cdot v_1| > (1- \alpha)|v|$}.
\label{hyp.rewrite2}
\end{equation}
Then there exist constants $a_0, \delta_2>0$ such that if $a\le a_0$ in \eqref{e1},
then
\begin{equation}
\label{wew55}
\int_{B_{\delta_2}} d\lt(Du,SO\lt(n\rt)A_1\rt) dx\leq C \ep.
\end{equation}
\end{a1}
\begin{proof}

By defining $\wt{A_i}:=A_1^{-T} A_i$ and $\ti{v}=A_1 v_1$ we find
hypothesis (\ref{hyp.rewrite2}) is preserved for the wells
$SO\lt(n\rt)\wt{A_1}\cup \dots SO\lt(n\rt) \wt{A_m}$ and so without
loss of generality we can assume $A_1=Id$.

Let $w$ denote the Lipschitz approximation of $u$ found in Proposition \ref{L9}.
Recall that in Lemma \ref{L12.5} we found a set $D\subset B_1$, with $L^n(D) \le C \ep$,
and such that $w$ is one-to-one and  $\det Dw(x) > \sigma$ in $B_{r_0}\setminus D$.

Note also that there exists a constant $C$ such that
\begin{equation}
\mbox{ if $x\not\in D$, then
$d(Dw(x)^{-1}, A_j^{-1}SO(n)) \le
C d(Dw(x), SO(n) A_j)$ for every $j$.}
\label{pf2.n1}\end{equation}
This is clear, because the fact that $w$ is Lipschitz
implies that $\{ Dw(x): x\not\in D\}$ is contained in
the compact set $\{ M : \det M>\sigma, \| M\| \le C(K) \}$.

Define also
\[
\YI := \{ x\in B_1 \setminus D \  : \ d(Dw(x), SO(n))  > d(Dw(x), K) \}.
\]
It suffices to prove that
\begin{equation}
L^n( \YI \cap B_{\delta_2}) \le C \ep
\label{LL2.sts}\end{equation}
for a suitable $\delta_2$; this readily implies \eqref{wew55}.

{\em Step 1}:
We first claim that for $L^{2n}$ a.e.  $(x,y)\in (B_{r_0}\setminus D)\times (B_{r_0}\setminus D)$
such that $\nu := \frac{y-x}{|y-x|}$ satisfies $|\nu  \cdot v_1| > 1-\alpha$,
we have the estimate
\begin{equation}
|x-y|\le |w(x) - w(y)| 
\ + \
\int_{[w(x), w(y)] }  \Theta  \ dH^1 \ - c\alpha H^1([w(x), w(y)]\cap w(\YI))
\label{pf2.claim1a}\end{equation}
where $\Theta$ is a nonnegative function, independent of
$x$ and $y$ and given explicitly below, such that
\begin{equation}
\int_{\beta_{\rho_0}} \Theta(\xi) \ d\xi \ \le C \ep,
\label{pf2.claim1b}\end{equation}
where recall ball $\beta_{\rho_0}$ is the large ball in the image
we obtain from Proposition \ref{L9} \em (v).\rm

To prove this, we use Lemma  \ref{Lpath}, which implies that for
a.e.  $(x,y)\in (B_{r_0}\setminus D)\times (B_{r_0}\setminus D)$
there is an injective Lipschitz path
$g:[0,1]\to \R^n$,  such that $g(0) = x, g(1) = y$, and $w(g(t)) \in [w(x), w(y)]$ for all $t\in [0,1]$.
We will write $\gamma(t) := w(g(t))$, so that $Dw(g(t)) g'(t) = \gamma'(t)$.
Then
\[
|y-x| = \nu^T (y-x)
= \nu^T \int_0^1 g'(t) dt
=  \int_0^1 \nu^T Dw(g(t))^{-1} \gamma'(t) dt.
\]
Let us temporarily write $M(t)$, or simply $M$, for $Dw(g(t))$.
Then applying the area formula to the right-hand side above
(since the image of $\gamma$ is the segment  $[w(x), w(y)]$), we deduce that
\[
|x-y|  \ = \
\int_{[w\lt(x\rt),w\lt(y\rt)]} \sum_{ \{ t\in [0,1]: \gamma(t)= \xi \} } \nu^T M(t)^{-1}  \frac{\gamma'(t)}{|\gamma'(t)|}  \ dH^1( \xi).
\]
It follows that
\begin{equation}
|x-y| = |w(x) - w(y)| + \int_{[w(x), w(y)] } \Theta_0(\xi) \ dH^1(\xi)
\label{inv.est1}\end{equation}
for
\begin{equation}
\Theta_0(\xi) =\left(\sum_{ \{ t\in [0,1]: \gamma(t)= \xi \} } \nu^T M(t)^{-1}  \frac{\gamma'(t)}{|\gamma'(t)|} \right) -1.
\label{Th0.def}\end{equation}
If $g(t)\in \YI$, then  there exists some $j\ge 2$ such that
$d(M(t),K) = d(M(t),SO(n)A_j)$. Recalling that $A_1=Id$, we infer from \eqref{hyp.rewrite2}  that
$|\nu^T A_j|^{-1} \le 1-\alpha$, so that
\begin{eqnarray*}
\nu^T M\lt(t\rt)^{-1}\frac {\gamma'}{|\gamma'|}
&\leq
& \lt|\nu^T M\lt(t\rt)^{-1}\rt| \\
&\leq&
|\nu^T A_j^{-1}| + d\lt(M\lt(t\rt)^{-1},A_j^{-1} SO\lt(n\rt)\rt) \\
&\overset{ (\ref{pf2.n1})}{\leq}&
1-\alpha + C d(M\lt(t\rt), K)
\quad\quad
\quad\quad
\quad\quad
\mbox{ if }g(t)\in \YI.
\end{eqnarray*}
Similarly,
\[  
\nu^T M\lt(t\rt)^{-1} \frac{\gamma'}{|\gamma'|}\le1+ d(M\lt(t\rt)^{-1}, SO(n) )\le1 + C d(M\lt(t\rt), K)\quad\quad
\quad\quad\text{  if
}g(t)\not\in (\YI\cup D).
\]  
For $g(t)\in D$, we claim that the fact that $\| M(t)\|\le C(K)$ implies that
\[  
\nu^T M\lt(t\rt)^{-1} \gamma' \le C |\det M(t)|^{-1}|\gamma'|
\]   
To see this, we recall the polar decomposition $M\lt(t\rt)=QS$, where $Q\in O(n)$ and $S = \sqrt {M^T M}$ is
symmetric and nonnegative. Then $|\nu^TM\lt(t\rt)^{-1}\gamma' | \le  |\nu^T S^{-1}| \, |\gamma'|
\le C(\max \{ \lambda_i^{-1} \})   |\gamma'|$, where $\{ \lambda_i\}$
are the eigenvalues of $S$. The fact that $M$ is bounded implies that
$\lambda_i^{-1} \ge C$ for all $i$, and it follows that $\max_i  \{ \lambda_i^{-1} \} \le
C \det S^{-1} = C|\det M|^{-1}$. This proves the claim.

Since $w$ is one-to-one on $B_{r_0}\setminus D$, the above estimates of $\nu^T M\lt(t\rt)^{-1}\frac{\gamma'}{|\gamma'|}$ imply
that
\[
\Theta_0(\zeta) \ \le - \alpha \chara_{\zeta \in w(\YI)} + \Theta(\zeta)
\]
where
\begin{equation}
\label{wei21}
\Theta(\zeta) :=
\chara_{\zeta \not\in w\lt(V\rt)} d(Dw(w^{-1}(\zeta)), K)
+c\chara_{\zeta\in w\lt(V\rt)} \sum_{w(z)=\zeta} |\det Dw(z)|^{-1}.
\end{equation}
We now see that \eqref{pf2.claim1a} follows from the above with \eqref{Th0.def}, \eqref{inv.est1}.
To verify \eqref{pf2.claim1b}, note that by a change of variables, Proposition \ref{L9} (v)(b)
(recall definition (\ref{xy43}))
\[
\int_{\beta_{\rho_0}} \chara_{\zeta \not\in w\lt(V\rt)} d(Dw(w^{-1}(\zeta)), K) \ d\zeta
=
\int \chara_{z \not\in D} d(Dw(z), K) \ |\det Dw(z)| dz \ \le C \int d(Dw,K) \ \le C\ep.
\]
And similarly,
\[
\int_{\beta_{\rho_0}} \chara_{\zeta\in w\lt(V\rt)} (\sum_{y\in w^{-1}\lt(\zeta\rt)} |\det Dw(y)|^{-1} ) \ d\zeta \ = \
\int \chara_{z\in D} \ dz = L^n(D) \le C \ep.
\]

{\em Step 2}.
By arguing exactly as in Step 1 of the proof of Lemma \ref{LL1} we can find points
$x_0,\ldots, x_n \in B_{1/2}\setminus D$, a number $0< \delta< 1/8$,
and an affine map $l_Q$ with $Dl_Q = Q\in SO(n)$, such that $|x_k|\ge 3/8$ for all
$k$, and all of
the following hold. First, if $x\in B_{2\delta}$ then
\begin{align*}
|\tau_k \cdot  v_1| \ge 1- \alpha
\quad\mbox{ for }\tau_k = \frac {x-x_k}{|x-x_k|},\ \ k=0,\ldots,n,\quad\mbox{ and }\\
x = \sum_{k=0}^n\lambda_k x_k\mbox{ with }\sum_{k=0}^n\lambda_k=1 \ \ \mbox{ and } \ \
\lambda_k>\delta \mbox{ for all }k.
\end{align*}
Second, $|w(x_k) - l_Q(x_k)| \le C \ep$. Third, \eqref{pf2.claim1a} holds for $(x,y) = (x_k, x)$, for $L^n$ a.e. $x\in B_{2\delta}\setminus D$.
And finally,
\begin{equation}
\label{wei18}
\int_{\beta_{\rho_0}} \Theta(\zeta) |\xi_k - \zeta|^{1-n} \ d\zeta  \le C \ep
\quad\quad\mbox{ for }\xi_k := w(x_k), k=0,\ldots, n.
\end{equation}

{\em Step 3}. We have defined $\xi_k := w(x_k)$ for $k=0,\ldots, n$.
We claim that for $\xi \in w(B_{2\delta}\setminus D)$
\begin{equation}
\sum_{k=0}^n  \int_{[\xi_k,\xi]} \chara_{ w(\YI) } \,dH^1  \le C \sum_{k=0}^n \int_{[\xi_k,\xi]} \Theta \ dH^1 \ + \ C \ep.
\label{pf2.step3}\end{equation}
We will write $e_k(\xi) := \int_{[\xi_k,\xi]} \Theta \ dH^1$.
Let $x=w^{-1}\lt(\xi\rt)$, since $w$ is injective on $B_{\delta}\backslash D$, point $x$ is well defined.
We first use
\eqref{pf2.claim1a} to see that
\begin{align*}
H^1( [\xi_k, w(x)] \cap w(\YI) \ )
&\le | w(x) - \xi_k| - |x-x_k| + e_k(\xi)\\
&\le
| l_Q(x) - w(x)| + |l_Q(x-x_k)| + |l_Q(x_k) - \xi_k| - |x-x_k| +e_k(\xi)\\
&=
| l_Q(x) - w(x)| + C \ep + e_k(\xi).
\end{align*}
To estimate $| l_Q(x) - w(x)|$ we argue as follows.  Since $l_Q$ is an isometry,
\begin{equation}
\label{weib12}
|l_Q( x) - w(x_k)|
\le |l_Q( x  - x_k)| + | w(x_k)-l_Q(x_k) |
\le |x_k-x| + C \ep.
\end{equation}
So we use \eqref{pf2.claim1a} again to find that
\begin{eqnarray}
\lt| l_Q(x) - w(x) -  (w(x_k) - w(x))\rt| &\overset{(\ref{weib12})}{\leq}& \lt|x_k-x\rt|+ C\ep\nn\\
&\overset{(\ref{pf2.claim1a})}{\leq}&|w(x)-w(x_k)| + e_k(\xi) + C \ep.\nn
\end{eqnarray}
It follows that
\[
\frac{ w(x) - w(x_k)}{|w(x)-w(x_k)| }\cdot ( l_Q(x) - w(x)) + |w(x_k)-w(x)|
\ \le \
|w(x)-w(x_k)| + e_k(\xi) + C \ep
\]
and hence that $\frac{w(x) - w(x_k)}{|w(x)-w(x_k)|}\cdot ( l_Q(x) - w(x)) \leq e_k(\xi)  + C \ep$
for $k=0,\ldots, n$.
Recall from Proposition \ref{L9} {\em (iv)} that there exists an affine map
$l_R$ such that
\begin{equation}
\mbox{
 $Dl_R=R\in SO(n)$ and $\| w - l_R\|_{L^\infty} \le C a^{\frac 1{n+1}}$}.
\label{recall.lR}\end{equation}
It follows that the convex hull of $\{ \frac{ w(x) - w(x_k)}{|w(x)-w(x_k)| } \}_{i=0}^n$
contains a ball $B_b$ of radius $b$ bounded away from zero, if  $a$ is small
enough. Thus Lemma \ref{AX1} implies that
$|  l_Q(x) - w(x)| \ \le C \sum_{k=0}^n  e_k(\xi)  + C \ep$,
and we have proved \eqref{pf2.step3}.

{\em Step 4}.  We use the notation
$
\Delta = w(D)$
and  we write $\beta_{\delta}: =  l_R(B_{\delta})$.
Note that $\beta_\delta$ is just a ball of radius $\delta$ (although not
necessarily centered at the origin), since
we are assuming that $A_1 = Id$. Note that from (\ref{recall.lR}) and (\ref{xy43}) we have
$\beta_{\delta}\backslash \Delta\subset w\lt(B_{2\delta}\backslash D\rt)$.

We next integrate \eqref{pf2.step3} over $\xi\in \beta_{\delta}\setminus \Delta$.
Both sides of the resulting inequality
contain terms of the form
$\int_{\beta_{\delta} \setminus \Delta} F_k(\xi) \ d\xi$, for $F_k$ of the form $F_k(\xi) = \int_{[\xi_k,\xi]} f(\eta) dH^1$.
Arguing exactly as in the proof of Step 3 of Lemma \ref{LL1} we find that
\begin{equation}
\label{wei6}
\int_{\beta_{\delta}\setminus \Delta} F_k(\xi) \ d\xi = \int f(\xi) G_k(\xi) |\xi-\xi_k|^{1-n} \ d\xi,\quad
\mbox{  where }
G_k(\xi_k+s\theta):=\int_s^\infty r^{n-1} \chara_{{\xi_k+r\theta\in \ \beta_\delta \setminus \Delta}} dr.
\end{equation}
Note in particular  that $G_k$ is bounded. We also claim that
\begin{equation}
\BI_k:=\{ \xi\in \beta_{\delta/2} \ : \ G_k(\xi) \leq c_0 \}\text{ is such that }
L^n\lt(\BI_k\rt)\le C_0 \ep
\label{pf2.step4}\end{equation}
for suitable constants $c_0,C_0$. To see this, fix
$\xi\in \beta_{\delta/2}$, and write $\xi$ in the form $\xi= l_R(x)$ with $|x| <\delta/2$.
Then for any $k$, since by definition $\xi_k = w(x_k)$,
\[
|\xi - \xi_k|
=
|l_R(x) - l_R(x_k) + l_R(x_k) - w(x_k)|
\overset{\eqref{recall.lR}}{\ge} |x-x_k| - C a^{\frac 1{n+1}}.
\]
It follows that if $a$ is small enough, then $|\xi - \xi_k| \ge 1/8$,  say, for $\xi\in \beta_{\delta/2}$.
Thus for $\theta\in S^{n-1}$ and $s>\frac{1}{8}$ such that $\xi_k=\xi+s\theta$
\begin{align*}
G_k\lt(\xi_k+s\theta\rt)
&\geq s^{n-1}\int_{s}^{\infty}
\chara_{\xi_k+r\theta\in \beta_{\delta} \setminus \Delta} dr
\geq 8^{1-n} L^1\lt(\lt\{r\in\lt[s,\infty\rt]:\xi_k+r\theta\in \beta_{\delta} \setminus \Delta \rt\}\rt)\\
&\ge 8^{1-n} \left[L^1\lt(\lt\{r\in\lt[s,\infty\rt]:\xi_k+r\theta\in \beta_{\delta} \rt\}\rt)
-  L^1\lt(\lt\{r\in\lt[s,\infty\rt]:\xi_k+r\theta\in  \Delta \cap \beta_{\delta} \rt\}\rt)\right].
\end{align*}
Any ray starting at a point in $\beta_{\delta/2}$
must travel a distance at least $\delta/2$ before leaving $\beta_\delta$, so the first term on the right-hand side above is greater than $\delta/2$ for $\xi_k+s\theta\in \beta_{\delta/2}$ . Thus \begin{eqnarray}
\label{wei11.32}
G_k\lt(\xi_k+s\theta\rt)
&{\geq}&
8^{1-n}
\lt(\delta/2-
L^1\lt(\lt\{r\in\lt[s,\infty\rt]:\xi_k+r\theta\in  \Delta \cap \beta_{\delta} \rt\}\rt)\rt)
\end{eqnarray}
for $\xi_k+s\theta\in \beta_{\delta/2}$.
Let $\Psi:=\lt\{\theta:  G_k\lt(\xi_k+s\theta\rt)\leq
8^{1-n} \frac{\delta}{4}
\ \ \mbox{ for some }\xi_k+s\theta\in \beta_{\delta/2} \rt\}$. Note that from (\ref{wei11.32})
\begin{equation}
\label{weib14}
L^1\lt(\lt\{r\in\lt[s,\infty\rt]:\xi_k+r\theta\in  \Delta \cap \beta_{\delta} \rt\}\rt)
\geq 8^{1-n}
\frac{\delta}{4}\text{ for all }\theta\in\Psi.
\end{equation}
We take $c_0=8^{1-n} \frac{\delta}{4}$, and we use the notation $l_{\theta}^z=\lt\{z+\lm \theta:\lm>0\rt\}$.
>From the definition (\ref{pf2.step4}) of $\BI_k$, we compute (using Lemma \ref{AX3})
\begin{eqnarray}
\label{wei22}
L^n\lt(\BI_k\rt)
&\leq&
\int_{\beta_{\delta/2} } \chara_{\BI_k}\lt(\xi\rt)\lt|\xi -\xi_k\rt|^{1-n} d\xi
\leq\int_{\theta\in\Psi}\int_{l_{\theta}^{\xi_k}\cap \beta_{\delta/2} } \chara_{\BI_k}\lt(\xi\rt) dH^1 \xi
dH^{n-1}\theta\nn\\
&\leq& H^{n-1}\lt(\Psi\rt) \nn \\
&\overset{(\ref{weib14})}{\leq} &
\frac{1}{c_0}
 \int_{\theta\in\Psi}\int_{l_{\theta}^{\xi_k}} \chara_{\beta_{\delta}\cap \Delta}\lt(\xi\rt) dH^1 \xi
dH^{n-1}\theta\nn\\
&=&
\frac{1}{c_0}\int \chara_{\Delta\cap \beta_{\delta} }\lt(\xi \rt)\lt|\xi -\xi_k\rt|^{1-n} d\xi
\leq
C L^n\lt( \Delta \rt)\leq  C \ep.\nn
\end{eqnarray}
Thus (\ref{pf2.step4}) is established.

{\em Step 5}. Defining $f\lt(\eta \rt)=\chara_{w\lt(\YI\rt)}\lt(\eta \rt)$ and
$F_k\lt(\xi\rt)=\int_{\lt[\xi_k,\xi\rt]} f\lt(\eta\rt) dH^1 \eta$ recall the definition of $\Theta$
from (\ref{wei21}), let $H_k\lt(\xi\rt):=\int_{\lt[\xi,\xi_k\rt]} \Theta\lt(\eta\rt) dH^1 \eta$ for
$\xi\in w\lt(B_{2\delta}\backslash D\rt)$, recall also that $\beta_{\delta}\backslash \Delta\subset
w\lt(B_{2\delta}\backslash D\rt)$
\begin{eqnarray}
\label{wzx1} C \ep+C\sum_{k=0}^n\int_{\beta_{\delta}\setminus \Delta } H_k\lt(\xi\rt) d\xi&\overset{(\ref{pf2.step3})}{\geq}&
\sum_{k=0}^n \int_{\beta_{\delta}\setminus \Delta } F_k\lt(\xi\rt) d\xi\nn\\
&\overset{(\ref{wei6})}{\geq}&
\sum_{k=0}^n\int_{\beta_{\delta} } f\lt(\xi\rt) G_k\lt(\xi\rt)\lt|\xi-\xi_k\rt|^{1-n} d\xi\nn\\
&\overset{(\ref{pf2.step4})}{\geq}&
c_0\int_{\beta_{\delta/2}} f\lt(\xi\rt) d\xi-C \ep.
\end{eqnarray}
Note that $\sum_{k=0}^n\int_{\beta_{\delta}\setminus \Delta } H_k\lt(\xi\rt) d\xi
\leq C \sum_{k=0}^n \int_{\beta_{\rho_0}} \Theta\lt(\xi\rt)\lt|\xi-\xi_k\rt|^{1-n} d\xi\overset{(\ref{wei18})}{\leq}
C\ep$,
so putting this together with (\ref{wei22}) we have
\begin{equation}
\label{wei2.87}
L^n(\beta_{\delta/2}\cap w(\YI))  \ = \ \int_{\beta_{\delta/2}} f\lt(\xi\rt) d\xi
\ \overset{(\ref{wzx1})}{\leq} \ C\ep.
\end{equation}
We remark also that \eqref{recall.lR} implies that
$w\lt(B_{{\delta}/{4}} \cap \YI \rt)\subset\beta_{\delta/2}\cap w(\YI)$ if $a$ is sufficiently small,
so that $L^n \lt(w\lt(B_{{\delta}/{4}} \cap \YI \rt) \rt) \le C \ep$.
Next recall that $\det Dw \ge \sigma$ in $\YI$, since $\YI \cap D = \emptyset $.
Now as $w$ is injective on $B_{\frac{\delta}{4}}\backslash D$,
the area formula implies that
\[
C \ep  \overset{\eqref{wei2.87}}{\ge} L^n\lt(w\lt(B_{\delta/4}\cap \YI\rt)\rt)
=
\lt|\int_{B_{\delta/4}\cap \YI} \det\lt(Dw\lt(x\rt)\rt) dx\rt|
\ \ge \ \sigma L^n(B_{\delta/4} \cap \YI),
\]
which is \eqref{LL2.sts}.
\end{proof}


\subsection{Sharp $L^{\infty}$ control on a large subset}

The methods used above yield the following result, which is valid for $m$ wells in $\R^n$
without any conditions on the wells.

\begin{a5}
\label{L20} Let $L=\cup_{i=1}^m SO\lt(n\rt)A_i$. Suppose  $u:B_1\to \R^n$  satisfies (\ref{e1}), and assume as
in \eqref{e3} that $A_1$ is the majority phase.
Let $r_1 = 1 - k_1a^{\frac 1{n+1}}$ be the constant found in
Proposition \ref{P.invert} and let $\ep=\vs^{\frac{1}{p}}$.

Then there exists $\OI\subset
B_{r_1}$ and some $R\in SO\lt(n\rt)A_1$ where
$L^n\lt(B_{r_1}\backslash \OI\rt)\leq C a^{1/p}$, and
\begin{equation}
\label{eq46} \|u-l_Q\|_{L^{\infty}\lt(\OI\rt)}\leq C \ep.
\end{equation}
\end{a5}

\begin{proof}
It suffices to prove the Proposition for all $a< a_0$, for some fixed $a_0>0$.

Let $\{ y_0,\ldots, y_n\}$ be the vertices of a regular simplex centered at $0$
with  $|y_k| = \frac 12 \ \forall k$.
Then using Corollary \ref{C.rsimp} (which is valid if $a_0$ is taken to be small enough)
we find points $x_k\in B_{1/8}(y_k)$
such that
\begin{equation}
\large\left| \ |u\lt(x_k\rt)-u\lt(x_l\rt)| \  -  \  |A_1(x_k - x_l)|  \ \large\right| \le  C  \ep
\quad\mbox{ for all }k\ne l
\label{xk1}\end{equation}
and
\begin{equation}
L^n(  \ \left\{ y\in B_{r_1} \ \: (x_k, y)\not \in \GI \right) \ ) \ \le C a^{1/p},
\label{xk2}\end{equation}
where $\GI$ is the set found in Proposition \ref{L10}.

Now it follows from \eqref{xk1} and Lemma \ref{L12} (see section \ref{SLA}) that there exists an affine
map $l_R$ with $R\in O(n)A_1$
such that
\[
| u(x_k) - l_R(x_k)| \le C \ep \quad\mbox{ for }k=0,\ldots,n
\]

Let
\[
\OI := \bigcap_{k=0}^n \{ y\in B_{r_1} \ : \ (x_k,y)\in \GI  \}
\]
It is clear from \eqref{xk2} that $L^n(B_{r_1}\setminus \OI) \le C a^{1/p}$. Next, note that
the definition of $\GI$ implies that
if $y\in \OI$, then
\begin{equation}
\large\left| |u\lt(y\rt)-u\lt(x_k\rt)| \  -  \  |A_1(y-x_k)| \large\right| \le  C \ep
\quad\quad\quad\mbox{ for }k=0,\ldots,n.
\label{OI.g1}\end{equation}
Then it  follows directly from the final
conclusion of Lemma \ref{L12} that $|u(y) - l_R(y)| \le C \ep$. This proves
\eqref{eq46}.

To complete the proof of Proposition \ref{L20} we only need to note that
by Proposition \ref{L9}, \it (iv) \rm there exists some $R\in SO\lt(n\rt)A_1$ and
affine function $l_R$ with $Dl_R=R$ such that
$\|w-l_R\|_{L^{\infty}\lt(\OI\rt)}\leq c a^{\frac{1}{n+1}}$. So in
particular
\begin{eqnarray}
\|l_{O}-l_R\|_{L^{\infty}\lt(\OI\rt)}&\leq&
\|w-l_{O}\|_{L^{\infty}\lt(\OI\rt)}+\|w-l_R\|_{L^{\infty}\lt(\OI\rt)}\nn\\
&\leq&a^{\frac{1}{n+1}}.\nn
\end{eqnarray}
Thus $O$ and $R$ must belong to the same connected component of $O\lt(n\rt)A_1$.
Therefore $O\in SO\lt(n\rt)A_1$.

\end{proof}

\section{Totally rank-$1$ connected wells}\label{sctwt}

Recall that we have shown that
an $m$-well Liouville Theorem holds for
$K = \cup_{i=1}^m SO(n)A_i$
satisfying the condition
\begin{equation}
\forall \ \ \mbox{  $i\in\lt\{1,2,\dots m\rt\}, \ \exists \  v_i\in S^{n-1}$
such that
$\lt|A_i v_i\rt|>\lt|A_j v_i\rt|\text{ for all }j\ne i$.}
\label{cond}\end{equation}
Given $K = \cup_{i=1}^n SO(n)A_i$,
we form a graph $\GI_K$ with vertices $v_1,v_2,\dots v_n$ where $\GI_K$
has edge $\lt(v_i,v_j\rt)$ if and only if there exists $A\in SO\lt(n\rt)A_i$ and $B\in SO\lt(n\rt)A_j$
with $\mathrm{rank}\lt(A-B\rt)=1$. We say that  $K$ is {\em totally rank-$1$ connected}
if $\GI_K$ forms a connected graph.
We will prove that, loosely speaking, \eqref{cond} is satisfied for
most totally rank-$1$ connected collections of $n$ wells in $\R^n$.


We say that a well $SO(n)A$ is  {\em positive} if $\det A>0$.
We will restrict our attention to positive wells, since we are interested in orientation-preserving
maps with nonvanishing determinant.
The map $SO(n) A \mapsto A^T A$ defines a bijection between the set of positive wells and the
set of positive definite symmetric matrices. This map is clearly well-defined, since
$\tilde A^T \tilde A = A^T A$ for any $\tilde A\in SO(n)A$, and it is invertible,
since $SO(n)A = SO(n) \sqrt{A^T A}$ when $SO(n)A$ is positive.

It is often convenient to describe properties of wells $SO(n)A$ in terms of the associated positive definite matrices
$A^T A$. An instance of this is the following well-known

\begin{a1}
\label{LL10}
Two positive wells $SO(n)A$ and $SO(n)B$ are rank-$1$ connected
if and only if there exist column
vectors
$p,q$, at least one of which is nonzero, such that $p\cdot q = 0$ and
\begin{equation}
A^TA - B^TB = p p^T - q q^T.
\label{r1c.qf}\end{equation}
\label{quad.r1c}\end{a1}

The matrix $p p^T - q q^T$ is uniquely determined by the two wells, so that the
wells determine the vectors $p,q$ up to multiplication by $-1$.
The degenerate cases $p=0$, $q=0$ are not excluded.
(Clearly if $p=q=0$ then $SO(n)A_1 = 
 SO(n) A_2$.)
We present a proof of Lemma \ref{quad.r1c}
at the end of this section, since we have not been able to find a good reference.

The main result of this section is the following

\begin{a2}
\label{T3}
Let $SO(n) A_i$ be positive wells for $i=1,\ldots, n$, and assume that there exists a set $\CI$ of the form
$\CI = \{ (i_k, j_k ) \}_{k=1}^{n-1}$
such that $SO(n)A_{i_k}$ is rank-$1$ connected to $SO(n)A_{j_k}$ for every $k$.
Assume moreover that
\begin{equation}
\mbox{ $\forall \ i,j\in \{1,\ldots, n\}$, \  $\ \exists w_0 = i, w_1,\ldots, w_\ell = j$
such that $(w_k,w_{k+1} )\in \CI$ or $(w_{k+1},w_k)\in \CI$.}
\label{CI}\end{equation}
For  $k=1,\ldots, n-1$, let $p_k, q_k$ be the vectors characterized
(up to a sign) by the conditions
\begin{equation}
\mbox{$p_k\cdot q_k=0$ and
$A_{i_k}^T A_{i_k} - A_{j_k}^T A_{j_k} =  p_k p_k^T - q_k q_k^T$},
\label{pk.qk}\end{equation}
and assume that
\begin{equation}
\mbox{$\{ p_1,\ldots, p_{n-1}, q_1,\ldots, q_{n-1} \}$
contains no linearly dependent subset of $n$ elements.}
\label{pk.qk1}\end{equation}
Then $K = \cup_{i=1}^n SO(n)A_i$ satisfies condition \eqref{cond}.
\label{T.generic}\end{a2}

\begin{remark}
\label{R1}
The assumptions about $\CI$ imply that $K$ is totally rank-$1$ connected,
since
$$
\lt\{\{ i_k, j_k\}:k=1,2,\dots n-1\rt\}
$$
are the edges of the graph $\GI_K$ that characterizes the rank-$1$
connectivity of $K$ and assumption \eqref{CI} implies that $\GI_K$ is connected.

Conversely, whenever $K$ is totally rank-$1$ connected, we can find a set
$\CI$ satisfying the above conditions. Indeed,
by an elementary result in graph theory, every connected graph with $n$ vertices has a connected
subgraph with the same vertices and only $(n-1)$ edges (this subgraph
is know as a \em spanning tree\rm). Thus given any totally rank-$1$ connected $K$,
we can select a spanning tree and use it to define $\CI$, by listing the
edges in some arbitrary order from $1$ to $n-1$,
and then orienting each edge by imposing an order on the associated vertices
(i.e.\, replacing the unordered pair $\{i_k,j_k\}$ by the ordered pair $(i_k,j_k)$ for example.)
\end{remark}

\begin{remark}

We claim that Theorem \ref{T.generic}  shows that \eqref{cond}
is satisfied  in
\[
\RI:= \{ K = \cup_{i=1}^n SO(n)A_i  \ :   K\mbox{ is totally rank-$1$ connected},  SO(n)A_i \mbox{ positive for all }i   \}
\]
except on a closed set of measure zero.
To see this, it suffices to argue that
for every set $\CI = \{ (i_k,j_k)\}_{k=1}^{n-1}$ satisfying \eqref{CI},
corresponding to a possible way of connecting the different wells,
the hypotheses of Theorem \ref{T.generic} are satisfied in
\begin{align*}
\RI_{\CI} &:= \{ \cup_{i=1}^n SO(n)A_i  :  SO(n)A_i \mbox{ rank-$1$ connected to }SO(n)A_j \mbox{ if }(i,j)\in \CI, \\
&\hspace{21em} SO(n)A_i \mbox{ positive for all } i\}
\end{align*}
away from a closed set of measure zero.
For simplicity we consider $\CI_0 = \{ (1,2), \ldots, (n-1, n)\}$,
corresponding to collections of wells such that $SO(n)A_i$
is rank-$1$ connected to $SO(n)A_{i+1}$ for $i=1,\ldots, n-1$. (The argument is
nearly identical for any other $\CI$.)
Consider the set
\[
\SI :=\{ (S, (p_1,q_1),\ldots, (p_{n-1}, q_{n-1})) : S\in M^{n\times n} \mbox{ is  symmetric, }
p_i, q_i\in \R^n, p_i\cdot q_i = 0 \mbox{  for all }i\}.
\]
Given $(S, (p_i,q_i)\, ) \in \SI$, we define symmetric matrices $S_1,
\ldots, S_n$ by
\[
S_1 := S, \quad\quad S_{i+1} := S_i +  p_i p_i^T - q_iq_i^T.
\]
Let  $\SI_+ := \{ (S, (p_i,q_i) ) \in \SI  : S_i\mbox{ as defined
above is positive definite for all } i \}$. Lemma \ref{quad.r1c}
implies that for $(S, (p_i,q_i) )\in \SI_+$, the collection $K =
\cup_{i=1}^n SO(n)\sqrt{S_i}$ belongs to $\RI_{\CI_0}$, and also
that every $K\in \RI_{\CI_0}$ arises in this fashion. Thus
$\RI_{\CI_0}$ can be parameterized by points in $\SI_+$, which is an
open subset of $\SI$.

The point is that  one can easily check that \eqref{pk.qk1} fails
only on a union of hypersurfaces in  $\SI$, which
is a closed set in $\SI$ of $H^{\dim \SI}$ measure $0$.

\end{remark}

\begin{proof}[Proof of Theorem  \ref{T.generic}]  Let $\HI_K$ be the graph with vertices
$v_1,v_2,\dots v_n$ where $\lt(v_i,v_j\rt)$ is an edge of $\HI_K$ if
and only $(i,j)$ or $(j,i)$ belong to $\CI$ \footnote{Note that in
general $\HI_K$ is only a subgraph of $\GI_K$}. As noted in Remark
\ref{R1}, $\HI_k$ is connected graph. We say an arbitrary connected
graph is a \em tree \rm if and only if it contains no \em loops\rm,
by this we mean it contains no non-trivial sequences of edges
$(v_{i_1},v_{i_2}), (v_{i_2},v_{i_3}), \dots (v_{i_{n-1}},v_{i_n})$
with $v_{i_1}=v_{i_n}$. It is well known (and easy to prove by
induction) that every tree with $n$ vertices has $(n-1)$ edges. Also
well known is that any connected graph contains a subgraph with the
same vertices that turns out to be a tree. From these two facts we
can conclude $\HI_K$ is itself a tree since it already has a minimal
possible number of edges.

{\em Step 1}.
We first claim that it suffices to show that for every
$(\sigma_1,\ldots, \sigma_{n-1}) \in \{\pm 1\}^{n-1}$, we can find a vector $v $
such that
\begin{equation}
\sigma_k ( |A_{i_k} v|^2 - |A_{j_k}v|^2) > 0
\quad\mbox{ for every }k= 1,m\ldots, n-1.
\label{r1c.suff}\end{equation}
We fix $i\in \{1,\ldots, n\}$. Let $\DI_1:=\lt\{j\in\lt\{1,2,\dots n\rt\}:(i,j)\in \CI\text{ or }(j,i)\in\CI\rt\}$.
For any $j\in\DI_1$ we chose the sign of $\sigma_k$ is the obvious way if
$(i,j)=(i_k,j_k)\in\CI$ then chose $\sigma_k=1$ and if $(i,j)=(j_k,i_k)\in\CI$ then $\sigma_k=-1$. Then
(\ref{r1c.suff}) gives us that $\lt|A_i v\rt|^2>\lt|A_j v\rt|^2$ for all $j\in\DI_1$.

Let $\DI_2:=\lt\{l\in\lt\{1,2,\dots n\rt\}:\text{For some
}p_l\in\DI_1, (l,p_l)\in\CI\text{ or }(p_l,l)\in\CI\rt\}$. For if
$l\in \DI_2$, if $(l,p_l)=(i_k,j_k)\in\CI$ chose $\sigma_k=1$ and if
$(p_l,l)=(j_k,i_k)\in\CI$ chose $\sigma_k=-1$. We then have $|A_i
v|^2>|A_{p_l} v|^2>|A_l v|^2$ for any $l\in \DI_2$. And we can
continue inductively defining $\DI_3,\DI_4,\dots$ choosing signs
such that (\ref{r1c.suff}) implies (\ref{cond}). We will never have
$i\in \DI_m$ (for any $m$) because the graph $\HI_K$ is a \em tree
\rm (and recall its (oriented) edges are given by $\CI$), for the
same reason if $m_2>m_1$ then $\DI_{m_2}\cap\DI_{m_1}=\emptyset$ and
the chain of inequalities we build will be consistent. The geometric
picture is that we start from a vertex on the graph $\HI_K$ and
expand outwards one edge at a time choosing signs $\sigma_k$ one at
the time.

{\em Step 2}.
We now fix an arbitrary $(\sigma_1,\ldots, \sigma_{n-1}) \in \{\pm 1\}^{n-1}$,
and we show that the system of inequalities \eqref{r1c.suff} admits a solution.
We  assume  for simplicity that $\sigma_k = 1$ for every $k$. This can be achieved
by replacing some pairs $(i,j)\in \CI$ by $(j,i)$; in fact the order
is arbitrary, and all our assumptions are preserved by this relabeling.
Then in view of the characterization of $p_k, q_k$,
\begin{eqnarray}
\lt(p_k\cdot v\rt)^2-\lt(q_k\cdot v\rt)^2&=&v^T p_k p_k^T v-v^T q_k q_k^T v=v^T A_{i_k} A_{i_k}^T v-v^T A_{j_k} A_{j_k}^T v\nn\\
&=&\lt|A_{i_k} v\rt|^2-\lt|A_{j_k} v\rt|^2\nn
\end{eqnarray}
so our task is  to find some
$v\in \R^n$ such that
\begin{equation}
(p_k\cdot v)^2 - (q_k\cdot v)^2 > 0
\label{r1c.suff2}\end{equation}
for every $k\in \{1,\ldots, n-1\}$.
To do this, fix a nonzero vector $v$  the subspace $\cup_{k=1}^{n-1} q_k^\perp$.
We then only need to check that $v\cdot p_k \ne 0$ for every $k$.
In fact, if $v\cdot p_k = 0$, then
$v\in (\mbox{span}(p_k, q_1, q_2, \ldots, q_{n-1}) )^\perp$.
Since $v\ne 0$, this would mean that $\{p_k, q_1, q_2, \ldots, q_{n-1} \}$
are linearly dependent, and in view of our assumptions this is impossible.
\end{proof}

We end this section by presenting a

\begin{proof}[Proof of Lemma \ref{quad.r1c} ]
First we claim that given arbitrary vectors $\ti{p},\ti{q}$ we can
find orthogonal $p,q$ such that $\ti{p}\ti{p}^T-\ti{q}\ti{q}^T=p
p^T-q q^T$. We write $M :=  \tilde p \tilde p^T -  \tilde q \tilde q
^T$. Note that $M$ is a symmetric matrix of rank at most $2$. This
is obvious, since $\tilde p \tilde p^T $ and $ \tilde q \tilde q ^T$
are both symmetric with rank $\le 1$. The claim is clear if
rank$(M)\le 1$, so we assume that rank$(M) = 2$. Then $\tilde p$ and
$\tilde q$ are linearly independent and in particular nonzero, so by
considering vectors orthogonal to $\tilde p$ and $\tilde q$
respectively one sees that $\min_{|v|=1} v^T M v < 0 < \max_{|v|=1}
v^TMv$. Thus $M$ has one positive and one negative eigenvalue, and
all other eigenvalues vanish (since rank$(M)=2$). The claim then
follows by diagonalizing $M$. So in fact it suffices to prove
$SO(n)A$ and $SO(n)B$ are rank-$1$ connected if and only if
(\ref{r1c.qf}) holds true for any $p,q\in\R^n$.

A second simplification comes from noting since $SO(n)B$ is
positive, $B$ is invertible, and after multiplying $A,B$ (on the
right)  and $p,q$ (on the left) by $B^{-1}$ we see that it suffices
to prove the lemma for $B = Id$.

Now assume that \eqref{r1c.qf} holds with $B = I$ and $p\cdot q = 0$.
The eigenvalues of $A^T A = I + pp^T - q q^T$
are $\lambda_1 := 1+|p|^2, \lambda_2 :=1 - |q|^2$, and $\lambda_3 = \ldots = \lambda_n = 1$.
Let $v_1,\ldots, v_n$ be an associated orthonormal basis of eigenvectors.
Since $A^TA$ is positive definite, it must be the case that $1-|q|^2 > 0$.
Let $A_0 = \sqrt {A^T A}$, so note that $A_0$ has an orthonormal basis of eigenvectors $v_1,\ldots, v_n$ and eigenvalues
$\mu_i  = \sqrt{\lambda_i}$ for $i=1,\ldots, n$. Note that $0< \mu_2 \le 1 \le \mu_1$, and that $\mu_1 < \mu_2$, since we have assumed that $A_0 \ne Id$.
To show that $A_0$ is rank-$1$ connected to some
$Q\in SO(n)$, it suffices to find $(n-1)$ orthonormal vectors
$w_1,\ldots, w_{n-1}$ such that $\{ A_0 w_i\}_{i=1}^{n-1}$ is also orthonormal, since then
we can take $Q$ to be the unique element of $SO(n)$ such that $Qw_i = A_0 w_i$ for
$i=1,\ldots, n-1$. Such a collection is provided by
\[
w_1 = \left(\frac{\mu_2^2-1}{\mu_2^2-\mu_1^2}\right)^{1/2}
 v_1 + \left(\frac{1- \mu_1^2}{\mu_2^2-\mu_1^2} \right)^{1/2} v_2,
\quad\quad\quad\quad w_i = v_{i+1}\quad\mbox{ for }i=2,\ldots, n-1.
\]

Now suppose that $A$ is rank-$1$ connected to
$SO(n)$, so that $A = Q +a b^T$ for some $Q\in SO(n)$ and nonzero column vectors
$a$, $b$. We can also assume that $1= |a|^2 = a^T a$; if not, replace $a$ by $\frac a{|a|}$ and $b$ by $|a|b$. Then
\[
A^TA - Q^TQ = Q^T ab^T + b a^T Q + b a^Ta b^T 
= \tilde ab^T + b \tilde a^T +  bb^T
\]
for $\tilde a  = Q^T a$. If we define $\tilde p = \tilde a +  b$ ,
it follows that $A^TA - Q^T Q = \tilde p \tilde p^T -  \tilde a \tilde a^T$.

And as we know we can find  orthogonal  vectors $p,q$ such that
$ \tilde p \tilde p^T -  \tilde a \tilde a ^T = p p^T - q q^T$, so that \eqref{r1c.qf} holds.
This finishes the proof of the lemma.
\end{proof}

\section{Appendix: Auxiliary Lemmas}\label{S:app}

\subsection{Truncation Lemma}

\begin{a1}
\label{L7}
Let $n,m\geq 1$, and suppose that $\Omega\subset \R^n$ is a bounded
Lipschitz domain. Suppose also that $f:\R^{m\times n}\to [0,\infty)$
is a function such that $f(v) \geq  c_1 |v|^{p} - c_2$ for some
$c_1>0, c_2\geq 0$, and $p\geq 1$. Then for any $q\in [1,\infty)$
there exists a constant $C$ 
such that, whenever  $u\in W^{1,q}(\Omega;\R^m)$ satisfies $f(Du)
\in W^{1,1}(\Omega;\R)$, then for every $\lambda>0$, there exists
$w\in W^{1,\infty}(\Omega;\R^n)$ such that
\begin{align*}
\it (i) \rm\quad\quad   &\| Dw\|_{L^\infty(\Omega)} \leq C\lambda, \\
\it (ii) \rm\quad\quad  &\| Du - Dw\|_{L^q(\Omega)}^q \leq
\frac C {\lambda^q} \int_{ \{x\in \Omega: Du(x)>\lambda \}} |Du|^q \ dx, \\
\it (iii) \rm\quad\quad &E := \{ x\in \Omega: u(x)\neq w(x)\}
\subset \{ x\in \Omega \ : \ \sup_r \dashint_{\Omega\cap B_r(x)} |Du| \,dy > \lambda \},\\
\it (iv) \rm\quad\quad  &|E | \leq
\frac C {\lambda^q} \int_{ \{x\in \Omega: Du(x)>\lambda \} }|Du|^q \ dx, \\
\it (v) \rm\quad\quad   &\mbox{ if }\ \ c_1\lambda^p-c_2>0, \mbox{ then  }\ \ Cap_1(E ) \leq \frac C{c_1\lambda^p-c_2 }\, \| f(Du)\|_{W^{1,1}(\Omega)}.
\end{align*}
Consequently, if $c_1\lambda^p-c_2>0$, then there exists an open set $E'$ with
smooth perimeter such that $E\subset E'$ and
$Per_{\Omega}(E')+\lt(L^n\lt(E'\rt)\rt)^{\frac{n-1}{n}} \leq \frac
C{c_1\lambda^p-c_2 }\, \| f(Du)\|_{W^{1,1}(\Omega)}$.
\end{a1}

We will apply the lemma with $f(Du) = d^p(Du,K)$ for some $p\geq 1$,
where $K$ is a compact subset of $\R^{m\times n}$.

Most of these conclusions are classical for $u\in
W^{1,q}_0(\Omega;\R^m)$ if $\Omega$ is smooth enough, and {\it (i) \rm,
\it (ii) \rm, \it (iv) \rm are proved in exactly the form stated above in
Proposition A.1, \cite{fmul}; hence we only sketch the proofs of
these points below. (These conclusions  do not require the
hypothesis $f(Du)\in W^{1,1}$.) The main point is {\em (v)}: control
over second derivatives of $u$ yields an estimate on the capacity of
the set $E = \{ x\in \Omega: u(x)\neq w(x)\}$.

If we assume $f(Du)\in W^{1,s}$ for some $s>1$, then by appealing to
slightly different results from the literature but otherwise leaving
the proof unchanged, we would obtain an estimate of $Cap_s(E)$. For
example, if $u\in W^{2,s}(\Omega)$ then (taking $f(Du) = |Du|$) we
would find that $Cap_s(E) \leq \frac C{\lambda^s} \| Du\|_{W^{1,s}}
^s$.

\begin{proof}
For any integrable function $v$ on any open subset $U\subset \R^n$,
we use the notation
\[
M_U(v)(x) := \sup_{r>0} \ \dashint_{B_r(x)\cap U} |v(y)| \  dy.
\]
For $u\in W^{1,1}(\Omega;\R^m)$ and $U\subset \Omega$ we will write
\[
R^\lambda(u; U ) := \{ x\in \Omega: M_U(Du)(x) \leq \lambda\}.
\]
We first assert that for any bounded Lipschitz domain $\Omega$,
there exists a constant $C$ such that for every $\lambda>0$ and
$u\in W^{1,1}(\Omega;\R^m)$,
\begin{equation}
|u(x) - u(y)| \leq C\lambda|x-y| \quad\quad \mbox{ for all }x,y\in
R^\lambda(u;\Omega). \label{L7.e1}\end{equation} This is well-known
if $\Omega = \R^n$ and is essentially proved in \cite{fmul} for
bounded Lipschitz domains; we recall the argument at the end of the
proof for the convenience of the reader. Once \eqref{L7.e1} is
known, standard extension theorems assert the existence of a
function $w:\Omega\to \R^m$ that satisfies the Lipschitz bound
\it (i) \rm and agrees with $u$ on $R^\lambda(u;\Omega)$, so that {\em
(iii)} holds. Then {\em (iv)} follows from {\em (iii)} by a covering
argument, and {\em (ii) } is a consequence of {\em (i), (iv)}.

To prove  {\em (v)}, we must estimate the $1$-capacity of
$\{x\in\Omega: M_\Omega(Du)(x) >\lambda\}$. To do this, note from
Jensen's inequality and the assumptions on $f$ that
\[
\dashint_{B_r(x)\cap \Omega} |Du| \ dy \leq \
\left(\dashint_{B_r(x)\cap \Omega} |Du|^p \ dy\right)^{1/p} \  \leq \
\left(\dashint_{B_r(x)\cap \Omega} \frac 1{c_1}( f(Du) + c_2) \
dy\right)^{1/p}.
\]
Thus  $\{x\in\Omega: M_\Omega(Du)(x) >\lambda\} \subset
\{x\in\Omega: M_\Omega(f(Du))(x) > (c_1\lambda^p - c_2)^+ \,\}$.
Hence {\em (v)} will follow once we check that
\begin{equation}
Cap_1 \left(\{x\in\Omega: M_\Omega(F)(x) >\mu \,\}\right) \leq \frac
C\mu \| F \|_{W^{1,1}(\Omega)} \label{L7.e2}\end{equation} for all
$F\in W^{1,1}(\Omega)$ and $\mu>0$. This is well-known, see
\cite{evans2} Section 4.8 for example, (and requires only  $\int
|DF|$ on the right-hand side) if $\Omega= \R^n$ and for example $F$
has compact support. To show that it remains valid in the present
circumstances, recall that any bounded, Lipschitz domain is an
extension domain (see for example \cite{stein}, Theorem 5 in Section
VI.3), so that there exists a function $\tilde F:\R^n\to \R$ with
support in a fixed compact set (independent of $F$), such that
\begin{equation}
\mbox{ $\tilde F = F$ on $\Omega$, and } \ \ \ \| \tilde F
\|_{W^{1,1}(\R^n)} \leq C \|  F \|_{W^{1,1}(\Omega)}.
\label{extension1}\end{equation} We may also take $\tilde F$ to be
nonnegative (since if this does not hold, we  may replace $\tilde F$
by $|\tilde F|$). Classical results mentioned above imply that
\[
Cap_1 \left(  \left\{ x\in \R^n : M_{\R^n}(\tilde F)(x) > \mu
\right\}\right) \leq \frac C \mu \int_{\R^n} |D\tilde F|,
\]
so in view of \eqref{extension1}, to prove \eqref{L7.e2} it suffices
to verify that
\begin{equation}
\mbox{ $M_\Omega(F)(x) \leq C M_{\R^n} (\tilde F)(x)$ for all $x\in
\Omega$}. \label{L7.e3}\end{equation} Fix a number $R >
\mathrm{diam}(\Omega)$, so that $\Omega\cap B_r(x) = \Omega$ if
$r\geq R$, for every $x\in \Omega$. Then for $x\in \Omega$,
\begin{align*}
M_\Omega(F)(x) &=
\sup_{0<r<R}\ \dashint_{\Omega\cap B_r(x)} F \, dy \\
&\leq \sup_{0<r<R}\
\frac 1{|\Omega\cap B_r(x)|} \int_{B_r(x)}\tilde F\, dy\\
&= \sup_{0<r<R}\ \frac {|B_r(x)|}{|\Omega\cap
B_r(x)|}\dashint_{B_r(x)} \tilde F \, dy \leq (\sup_{0<r<R}\frac
{|B_r(x)|}{|\Omega\cap B_r(x)|}) M_{\R^n}(\tilde F)(x).
\end{align*}

And the fact that $\Omega$ is Lipschitz implies that
$\sup_{0<r<R}\frac {|B_r(x)|}{|\Omega\cap B_r(x)|}  <\infty$; if
this were false, we could find a sequence of balls $B_k =
B_{r_k}(x_k)$, with $r_k$ necessarily tending to zero, with $x_k\in
\Omega$, and such that the density ratios $\frac{|\Omega\cap B_k|}
{|B_k|}$ tend to zero, and in a bounded  Lipschitz domain this is
easily seen to be impossible. Thus we have proved \eqref{L7.e2}, and
hence conclusion {\em (v)} as well.

To prove the final assertion about the existence of the set $E'$,
note that by the definition of capacity, there exists a function
$h\in C^{\infty}_c\lt(\R^n\rt)$ such that
\[
E \subset \mathrm{int}\lt\{x:h\lt(x\rt)\geq 1\rt\} \mbox{ and }
\int_{\R^n} \lt|Dh\rt|\leq 2 \,Cap_1(E).
\]

By the coarea formula
$$
\int_{0}^{1} H^{n-1}\lt(h^{-1}\lt(t\rt)\rt) dL^1
t\leq \int_{\R^n} |Dh| . 
$$
Thus we must be able to find $t_0\in\lt(1/2,1 \rt)$ with the
property that $H^{n-1}\lt(h^{-1}\lt(t_0\rt)\rt)\leq 4\, Cap_1(E)$.
We take $E' = \{ x\in \Omega: h(x) > t_0 \}$, so that the perimeter
estimate is satisfied. As in the proof of Lemma \ref{L3}, we can assume
$t_0$ is one of the a.e.\ numbers in $(\frac{1}{2},1)$ such that by Sard's theorem, $E'$ has
smooth boundary. And by Chebyshev and Sobolev inequalities, we
know that
\[
|E'|^{\frac {n-1}n} \leq C \| h \|_{L^{n/(n-1)}} \leq C \|Dh\|_{L^1}.
\]

Finally, we sketch the proof of \eqref{L7.e1}. If $\Omega$ is the
unit cube $Q$, then as noted by \cite{fmul}, one can deduce
\eqref{L7.e1} by minor modifications of classical arguments, as
expounded for example in \cite{evans2} chapter 6. Next, suppose that
$\Omega$ is a standard Lipschitz domain, or in other words, the
image of the unit cube under a map of the form $x = (x', x_n)\mapsto
\phi(x) = (x', q(x')) $, for $q:\R^{n-1}\to \R$ Lipschitz, note that
$\phi$ is a biLipschitz mapping. Then
given any $u\in W^{1,1}(\Omega;\R^m)$, we define $\tilde u:Q\to
\R^m$ by $\tilde u = u \circ \phi$. It is straightforward to check
that
\[
M_Q(D\tilde u)(x) \leq C M_\Omega(Du)(\phi(x)),
\]
and hence that \eqref{L7.e1} in this case follows from applying the
previous case to $\tilde u$. Finally, we note as in \cite{fmul} that
a bounded Lipschitz domain $\Omega$ can always be written as a
finite union $\Omega = \cup_{i=1}^k \Omega_i$, where each $\Omega_i$
is (up to a change of variables) a standard Lipschitz domain, so
that \eqref{L7.e1} holds for each $\Omega_i$. This can be done in
such a way that there exists some $r_1>0$ with the property that
there exists some $r_1>0$ such that,  for any $x,y\in \Omega$ such
that $|x-y|<r_1$, there exists some $i$ such that $\Omega_i$
contains both $x$ and $y$. Since $\Omega$ is bounded, it clearly
suffices to prove \eqref{L7.e1} for pairs $x,y$ such $|x-y|<r_1$, so
we need only show that for every $i= 1,\ldots,k$, there exists some
$C$ such that if $x,y\in \Omega_i\cap R^\lambda(u;\Omega)$, then
$|u(x)-u(y)|\leq C\lambda |x-y|$.

To do this, we fix some $i$ and argue as in the proof of
\eqref{L7.e3} above to find that
\[
M_{\Omega_i}(Du)(x) \leq C M_\Omega(Du)(x). \quad\quad\quad\mbox{ for
all $x\in \Omega_i$ and $u\in W^{1,1}(\Omega;\R^m)$.}
\]
Thus if $x,y\in \Omega_i\cap R^\lambda(u;\Omega)$ then $x,y\in
R^{C\lambda}(u;\Omega_i)$, and so the estimate $|u(x) - u(y)| \leq C
\lambda |x-y|$ follows from the case of standard Lipschitz domains.

\end{proof}

\subsection{Paths in the inverse of segments}

\begin{a1}
\label{Lpath}
Let $w:B_1\subset \R^n\rightarrow \R^n$ be a Lipschitz function.
Given a convex open set $\Lambda\subset \R^n$ such that
$\Lambda \subset w(B_1)\setminus w(\partial B_1)$,
let
\begin{equation}
\deg(w, B_1, \xi) = d_0 \ne 0\quad\mbox{ for a.e.  }\xi\in \Lambda
\label{piis.deg}\end{equation}
Then for $L^{2n}$ a.e. \ $(\eta,\zeta)\in\Lambda \times\Lambda $,
\begin{align}
&\exists\ \  b>0 \mbox{ and  an injective Lipschitz function }g:\lt[0,b\rt]\rightarrow B_1\mbox{ such  that }
\nn\\
&\hspace{7em}
w(g\lt(0\rt))=\eta, w(g\lt(b\rt))=\zeta,\mbox{ and }w\lt(g\lt(t\rt)\rt) \in \lt[\eta,\zeta\rt]
\quad\forall t\in [0,b].
\label{piis}\end{align}
\end{a1}

We will employ the framework of geometric measure theory, so that we
work with 
integral $k$-currents. One can think of such a current as a
$k$-submanifold of a Euclidean space that is described by specifying how it acts
(via integration) on $k$-forms. We will write
$\int_T\phi$ to indicate the action  of a current $T$ on a form $\phi$.
We will appeal to a number of classical facts
about slicing of currents. The basic reference for this material is \cite{fed} Chapter 4.3, and a more accessible discussion,
albeit without complete  proofs, can be found in \cite{gia}  section 2.5 of Chapter 2.

\begin{proof}

{\em Step 1}.
We will write $W(x) = (x, w(x))\in \R^n\times \R^n$ for $x\in B_1$, and $p_2((x,\xi)) =\xi\in \R^n$
for $(x,\xi)\in \R^n\times \R^n$. Note that $w = p_2\circ W$.

We write $G_w$ to denote  the (current associated with the) graph of $w$,
defined by
\[
\int_{G_w}\phi  \  
 := \ \int_{B_1} W^\# \phi
\]
for an $n$-form $\phi$ in $\R^n\times \R^n$, where $W^\#$ denotes the pullback via $W$.
(One can see $G_w$ as an example of a {Cartesian current}, and  an explicit expression for $G_w$ can be found on page 230, \cite{gia}.)
The boundary $\partial G_w$ of $G_w$ is defined by
$\int_{\partial G_w}\phi := \int_{G_w} d\phi$, and then the definition of $G_w$ implies that
$\int_{\partial G_w}\phi  = \int_{\partial B_1} W^\# \phi$.
These formulas  imply that
\[ 
\mbox{Spt}\,  G_w \ = \
\{ (x, w(x)) : x\in \bar B_1\},
\quad\quad\quad\quad
\mbox{Spt}\, \partial G_w \ = \
\{ (x, w(x)) : x\in \partial B_1\}.
\] 
We are using the fact that $w$ is Lipschitz, so that $\{ (x,w(x)) : x\in \mbox{ compact set }S\}$
is closed.

{\em Step 2}.
For $\nu\in S^{n-1}$ we define the functions
\[
q_\nu(\xi) := \xi - (\xi\cdot \nu) \nu = \mbox{ orthogonal projection onto }\nu^\perp\subset \R^n,\quad\quad\quad
Q_\nu  := q_\nu \circ p_2.
\]
We will write $\xi'$ to denote a generic point in $\mbox{Image}(q_\nu) =  \nu^\perp$.
We will need some classical results about slices of integral currents.
Recall that $\langle G_w, Q_\nu, \xi'\rangle$
denotes the slice of $G_w$ by $Q_\nu^{-1}(\xi')$, which for $H^{n-1}$ a.e. $\xi'\in \nu^\perp$
is a integral $1$-current satisfying
\[
\mbox{Spt}\, \langle G_w, Q_\nu, \xi'\rangle\subset \mbox{Spt} \,G_w\cap Q_\nu^{-1}(\xi'),
\quad\quad\quad 
\mbox{Spt} \, \partial  \langle G_w, Q_\nu, \xi'\rangle\subset \mbox{Spt}\,  \partial G_w\cap Q_\nu^{-1}(\xi')
\]
(see \cite{fed} 4.3.8 (2) for the first inclusion, 4.3.1 p437 together with 4.3.8 (2) for the second inclusion, alternatively
Section 2.5 \cite{gia} for a more readable presentation).
The fact that a.e.\ slice $\langle G_w, Q_\nu, \xi'\rangle$ is an integral $1$-current implies
(see \cite{fed} 4.2.25) that
we can write
\begin{equation}
\langle G_w, Q_\nu, \xi'\rangle = \sum_i R_i
\quad\quad\mbox{ for every $\nu$ and }H^{n-1} \ a.e.\  \xi'
\label{decomp}\end{equation}
where each $R_i = R_i(\nu, \xi')$ is the image of an injective Lipschitz map $\gamma_i: I_i\subset \R\to
\R^n\times \R^n$, and $I_i = (a_i, b_i)\subset \R$ is a bounded interval.
That is, $\int_{R_i} \phi = \int_{\gamma_i(I_i)} \phi =  \int_{I_i} \gamma_i^\# \phi$
for every $1$-form $\phi$ in $\R^n\times \R^n$.
The decomposition \eqref{decomp} is such that
\begin{equation}
\mbox{Spt}\, R_i \subset \mbox{Spt}\, \langle G_w, Q_\nu, \xi'\rangle \subset   \mbox{Spt}\, G_w\cap Q_\nu^{-1}(\xi')
\label{spt1}\end{equation}
\begin{equation}
\mbox{Spt}\, \partial R_i \subset \mbox{Spt}\, \partial \langle G_w, Q_\nu, \xi'\rangle
 \subset   \mbox{Spt}\, \partial G_w\cap Q_\nu^{-1}(\xi')
\label{spt2}\end{equation}
for every $i$.
It follows from \eqref{spt1}  that
each $\gamma_i$ has the form $\gamma_i(t) = (X_i(t), w(X_i(t)))$
for some Lipschitz path $X_i:I_i\to B_1$ such that $w(X_i(t)) \in q^{-1}_{\nu}(\xi')$
for every $t$.

The $0$-current $\partial R_i$ appearing in \eqref{spt2} is defined by  $\partial R_i(\phi) = \phi( \gamma_i(b_i^-)) - \phi(\gamma_i(a_i^+))$ for every function smooth
function $\phi$ on $\R^n\times \R^n$, so \eqref{spt2} asserts that if $\gamma_i(a_i^+) \ne \gamma_i(b_i^-)$ --- that is, if
$\partial R_i\ne 0$ --- then  $\gamma_i(a_i^+), \gamma_i(b_i^-)
\in \{ (x,w(x)) \ : \ x\in \partial B_1 \}$. In particular
\begin{equation}
\mbox{ if }\partial R_i\ne 0\mbox { then  }w(X_i(a_i^+)), w(X_i(b_i^-))\not\in \Lambda.
\label{spt3}\end{equation}

{\em Step 3}.
It $T$ is any $k$-current in $\R^n\times \R^n$,
we define $p_{2\#}G_w$ to be the $k$-current in $\R^n = \mbox{Image}(p_2)$ characterized by
$\int_{p_{2\#}T} \phi = \int_{T}p_2^\#\phi$, and we write $T\rest\Lambda$ to denote the restriction of
$T$ to $\Lambda$.
We claim that
\begin{equation}
(p_{2\#} \langle G_w, Q_\nu, \xi'\rangle ) \rest \Lambda=
 d_0 \langle \Lambda , q_\nu, \xi'\rangle
\label{piis.c3}\end{equation}
for every $\nu$ and $H^{n-1}$ a.e. $\xi'$, for $d_0$ as in \eqref{piis.deg}.
It follows from basic properties of slicing that the current on the right-hand side is just
the line segment $\Lambda \cap q_\nu^{-1}(\xi')$, with orientation and (nonzero)
multiplicity.

Since $Q_\nu  = q_\nu \circ p_2$,
\[
p_{2\#} \langle G_w, Q_\nu, \xi'\rangle =
p_{2\#} \langle G_w, q_\nu\circ p_2, \xi'\rangle =
 \langle p_{2\#}G_w, q_\nu, \xi'\rangle
\]
for a.e. $\xi'$, see \cite{fed} 4.3.2(7) for the last identity.
It follows that
\[
(p_{2\#} \langle G_w, Q_\nu, \xi'\rangle ) \rest \Lambda=  \langle (p_{2\#}G_w)\rest \Lambda, q_\nu, \xi'\rangle .
\]
Thus to prove \eqref{piis.c3}, it suffices to verify that
$ (p_{2\#}G_w)\rest \Lambda = d_0\Lambda$.
To prove this, we first note from the definitions that
\[
\int_{p_{2\#}G_w} \phi = \int_{B_1} W^\# p_2^\# \phi = \int_{B_1} (p_2\circ W)^\#\phi = \int_{B_1} w^\#\phi.
\]
In particular, if we write $\phi = \phi(\xi) d\xi$, where $d\xi$ denotes the standard volume form on $\R^n$,
then $w^\#\phi = \phi(w(x)) \det Dw(x) \,dx$, and so
the change of variables degree formula implies that
\[
\int_{p_{2\#}G_w} \phi\
= \int_{B_1} \phi(w(x)) \det Dw(x) \ dx\\
=
 \ \int_{\R^n} \phi(\xi) \deg( w, B_1, \xi) \, d\xi.
 \]
We conclude from \eqref{piis.deg} that
$
\int_{p_{2\#}G_w}  \phi\  = \  d_0 \int_{\Lambda} \phi$
if
$\mbox{Spt}\,\phi\subset \Lambda$. This says exactly that
$(p_{2\#}G_w)\rest \Lambda = d_0 \Lambda$, which is what we needed to prove.

{\em Step 4}.
We next claim that for every $\nu$, for $H^{n-1}$ a.e. $\xi'\in \nu^\perp$ and
every $i$ in the decomposition \eqref{decomp},
\begin{equation}
(p_{2\#} R_i) \rest \Lambda = d_i \langle \Lambda, q_\nu,
\xi'\rangle\quad\mbox{ for some }d_i\in \Zb.
\label{piis.c4}\end{equation} To see this, let us write $\Xi_i(t) =
p_{2}\circ \gamma_i(t) = w(X_i(t))$. Then it follows from the
definitions that $\int_{p_{2\#} R_i} \phi = \int_{I_i} \Xi^\#_i
\phi$. In view of properties of $\Xi_i = w\circ X_i$ recorded in
Step 2, this implies that $p_{2\#}R_i$ is supported in the line
segment $q^{-1}_\nu(\xi')$, and moreover \eqref{spt3} implies that
$\partial (p_{2\#}R_i) = 0$ in $\Lambda\cap q^{-1}_\nu(\xi')$. Then
\eqref{piis.c4} follows from the Constancy Theorem, see for example
\cite{fed} 4.1.7. (One can also deduce \eqref{piis.c4} by elementary
arguments from the fact that $\int_{p_{2\#} R_i} \phi = \int_{I_i}
\Xi^\#_i \phi$, together with the properties of $\Xi_i$ used above.)

{\em Step 5}.
It follows from \eqref{decomp} that for every $\nu$ and a.e. $\xi'\in\nu^\perp$,
$p_{2\#}\langle G_w, Q_\nu, \xi'\rangle = \sum_i p_{2\#} R_i$,
In view of Steps 3 and 4, this implies that the integer $d_i$ in \eqref{piis.c4} is
nonzero for at least one $i$.
Then the fact that $p_{2\#}R_i = d_i \langle \Lambda, q_\nu, \xi'\rangle$
implies that the for the corresponding curve $X_i$, the image of
$w\circ X_i$ covers $\Lambda\cap q^{-1}(\xi')$ and is contained in $q^{-1}_\nu(\xi')$.
So for any two points $\eta, \zeta$ in $\Lambda\cap q^{-1}(\xi')$,
we can find a path $g:[0,b]\to B_1$ satisfying \eqref{piis} by defining $g$
to be a reparametrization of the restriction of $X_i$ to a suitable subinterval of $I_i$.

{\em Step 6}.
Let $\BI:= \{ (\eta,\zeta)\in \Lambda\times\Lambda \ : \ \eqref{piis} \mbox{ does not hold}. \}$
Our goal is to show that $L^{2n}(\BI) = 0$. Note that from Step 5
\begin{equation}
H^2\left(\left\{ (\eta, \zeta)\in \BI : \eta,\zeta\mbox{ both belong to }\Lambda\cap q_\nu^{-1}(\xi') \right\} \right) = 0
\label{piis1}\end{equation}
for every $\nu\in S^{n-1}$ and $H^{n-1}$ a.e. $\xi'\in \nu^\perp$.

Let $f=\chara_{\BI}$. By Fubini's theorem
\begin{eqnarray}
\int_{\R^n}\int_{\R^n} f\lt(x,y\rt) dy dx&=&\int_{\R^n}\int_{\eta\in S^{n-1}}\int_{t>0}
f\lt(x,x+t\eta\rt) t^{n-1} dt dH^{n-1} \eta dx\nn\\
&=&\int_{\eta\in S^{n-1}}\int_{t>0}\int_{\R^n} f\lt(x,x+t\eta\rt) t^{n-1} dx dH^{n-1} \eta dt\nn\\
&=&\int_{\eta\in S^{n-1}}\int_{y\in \eta^{\perp}}\int_{s>0}\int_{t>0} f\lt(y+s\eta,y+\lt(t+s\rt)\eta\rt) ds dt
dH^{n-1} y dH^{n-1}\eta.\nn
\end{eqnarray}
>From (\ref{piis1}) for any $\eta\in S^{n-1}$, and $H^{n-1}$ a.e.\ $y\in \eta^{\perp}$
$\int_{s>0}\int_{r>0} f\lt(y+s\eta,y+\lt(t+s\rt)\eta\rt) t^{n-1} ds dt=0$ and thus
we have shown $L^n\lt(\BI\rt)=0$.
\end{proof}

\subsection{A linear algebra lemma}\label{SLA}

\begin{a1}

\label{L12} Suppose that $A$ is an invertible $n\times n$ matrix, and that
$z_0,z_1,\dots
z_n\in B_1\lt(0\rt)\subset \R^n$
and
$\zeta_0,\zeta_1,\ldots \zeta_n\in\R^n$
are points
such that
\begin{equation}
B_b\lt(y\rt)\subset\mathrm{conv}\lt(z_0,z_1,\dots z_n\rt)\text{ for
some }b>0, y\in B_1
\label{L12.h1}\end{equation}
and
\begin{equation}
\label{a6}
\lt|\lt|\zeta_i-\zeta_j\rt|-\lt|A(z_i-z_j)\rt|\rt|\leq \ep\text{ for all }
i\not=j\in\lt\{0,1,\dots n\rt\}.
\end{equation}

Then there exists an affine function $l_O$ with $Dl_O=O\in
O\lt(n\rt)A$ and constant $C=C\lt(b,n,A\rt)$ such that
\begin{equation}
\label{e10} \lt|\zeta_i-l_O\lt(z_i\rt)\rt|\leq C\ep\text{ for
all } i\in\lt\{0,1,\dots n\rt\}.
\end{equation}
Furthermore, if $z\in B_1$ and $\zeta\in \R^n$ are any other points
such that
\begin{equation}
\label{a6.1}
\lt|\lt|\zeta_i -\zeta \rt|-\lt|A(z_i-z)\rt|\rt|\leq \ep\text{ for all }
i \in\lt\{0,1,\dots n\rt\}
\end{equation}
then $|\zeta - l_O(z)| \le C \ep$ for the same $l_O$  as in \eqref{e10}, and with $C = C(n,b,A)$.
\end{a1}

\begin{proof}[ Proof of Lemma \ref{L12}. ]

By a translation we can assume that $z_0=\zeta_0=0$.
We can also assume that $A$ is the identity; if not,
simply replace each $z_i$ by $\tilde z_i = A z_i$ and drop
the tildes, so that \eqref{a6} becomes
\begin{equation}
\label{a6.2}
\lt|\lt|\zeta_i-\zeta_j\rt|-\lt|z_i-z_j\rt|\rt|\leq C\ep\text{ for all }
i\not=j\in\lt\{0,1,\dots n\rt\}.
\end{equation}
After these changes,
$|z_i|, |\zeta_i| \le C$ for all $i$.

We define $l_{\wt{O}}:\R^n\to \R^n$ to be the unique linear map
satisfying $l_{\wt{O}}\lt(z_i\rt)=\zeta_i$ for $i=1,2,\dots n$. We
will identify $l_{\wt{O}}$ with the matrix ${\wt{O}} = Dl_{\wt{O}}$.
It follows from \eqref{L12.h1} that  $\lt\{z_1,z_2,\dots z_n\rt\}$
are linear independent, and hence that $\wt{O}$ is well defined.

\em Step 1. \rm We first show that
\begin{equation}
\label{bc4}
\lt|\wt{O}\lt(z_i\rt)\cdot \wt{O}\lt(z_j\rt)-z_i\cdot z_j\rt|\leq C\ep.
\end{equation}
Toward this goal, note that since
$\lt|\ \lt|\zeta_i-\zeta_j\rt| + \lt|z_i-z_j\rt| \ \rt| \le C$
for all $i,j$,
\begin{eqnarray}
\label{xc40}
\lt|\lt|\zeta_i-\zeta_j\rt|^2-\lt|z_i-z_j\rt|^2\rt|
&\leq&c\lt|\lt|\zeta_i-\zeta_j\rt|-\lt|z_i-z_j\rt|\rt|\nn\\
&\overset{\lt(\ref{a6}\rt)}{\leq}& C\ep.
\end{eqnarray}
As a result,
\[
2\lt|\zeta_i\cdot \zeta_j-z_i\cdot z_j\rt|
\overset{(\ref{xc40})}{\leq}
\lt|\ |\zeta_i|^2+|\zeta_j|^2-|z_i|^2-|z_j|^2\ \rt|+C\ep.
\]
However, since $z_0 = \zeta_0 =0$, the $j=0$ case of \eqref{a6.2} implies that
$\lt|\ |\zeta_i|^2-|z_i|^2\ \rt| \le C\ep$,
and similarly
$\lt|\ |\zeta_j|^2-|z_j|^2\ \rt| \le C\ep$,
so
(\ref{bc4}) follows from the above.

\em Step 2. \rm We next claim that for any $v\in S^{n-1}$ there exist
$\gamma_1,\gamma_2,\dots \gamma_n$ with $\lt|\gamma_i\rt|\leq
\frac{2}{b}$ for each $i=1,2\dots n$ such that $\sum_{i=1}^{n}
\gamma_i z_i=v$.

\em Proof of Claim.  \rm Note that
\begin{align*}
B_b\lt(-y\rt)\cup B_b\lt(y\rt)
&\subset \mathrm{conv}\lt(z_0,z_1\dots
z_n\rt)\cup \mathrm{conv}\lt(z_0,-z_1\dots
-z_n\rt)\\
&\subset \mathrm{conv}\lt(z_1,\dots z_n,-z_1,\dots -z_n\rt)
\end{align*}
which implies $B_b\subset\mathrm{conv}\lt(z_1,\dots z_n,-z_1,\dots -z_n\rt)$.
So there exist positive $\beta_0,\beta_1,\dots \beta_{2n}$ such that
$\sum_{i=0}^{2n} \beta_i=1$ and $\sum_{i=1}^n \lt(\beta_i-\beta_{i+n}\rt)z_i = v b.$
Since $\frac{\lt|\beta_i-\beta_{i+n}\rt|}{b}\leq \frac{2}{b}$ this
completes  Step 2.

\em Step 3. \rm Let $\lt\{e_1,e_2,\dots e_n\rt\}$ be an orthonormal
basis of $\R^n$. We claim that
\begin{equation}
\label{bc6}
\lt|\wt{O}\lt(e_i\rt)\cdot \wt{O}\lt(e_j\rt)-\delta_{ij}\rt|\leq C\ep\text{ for any }i,j\in\lt\{1,2,\dots n\rt\}.
\end{equation}

\em Proof of Claim. \rm By Step 2 we can find coefficients
$\alpha^i_j\in \R$ such that $\sum_{j=1}^{n}  \alpha_j^i z_j=e_i$
and $\lt|\alpha^i_j\rt|\leq \frac{2}{b}$ for $i,j\in\lt\{1,2,\dots
n\rt\}$.  Note
\begin{equation}
\label{xc55} \sum_{k,l=1}^{n} \alpha^{i}_k\alpha^{j}_{l} z_k\cdot
z_l=\delta_{ij}\text{ for any }i,j\in\lt\{1,2,\dots n\rt\}.
\end{equation}
Now
\begin{eqnarray}
\lt|\wt{O}\lt(e_i\rt)\cdot
\wt{O}\lt(e_j\rt)-\delta_{ij}\rt|&\overset{(\ref{xc55})}{=}&\lt|\wt{O}\lt(e_i\rt)\cdot
\wt{O}\lt(e_j\rt)
-\sum_{k,l=1}^{n} \alpha^{i}_k\alpha^{j}_l z_k\cdot z_l\rt|\nn\\
&\leq&\lt|\sum_{k,l=1}^n  \alpha^i_k\alpha^j_l\lt(\wt{O}\lt(z_k\rt)\cdot\wt{O}\lt(z_l\rt)-z_k\cdot z_l\rt)\rt|\nn\\
&\leq&\sum_{k,l=1}^n   \lt|\alpha^i_k\alpha^j_l\rt|\lt|\wt{O}\lt(z_k\rt)\cdot\wt{O}\lt(z_l\rt)-z_k\cdot z_l\rt|\nn\\
&\overset{(\ref{bc4})}{\leq}& C\ep.\nn
\end{eqnarray}
Thus (\ref{bc6}) is established.

\em Step 4. \rm
We now define $\lt\{\xi_1,\xi_2,\dots \xi_n\rt\}$
to be the orthonormal basis of $\R^n$ obtained via a
Gram-Schmidt orthognalisation of the set of vectors
$\lt\{\wt{O}\lt(e_1\rt),\wt{O}\lt(e_2\rt),\dots
\wt{O}\lt(e_n\rt)\rt\}$.
Then an easy induction argument using \eqref{bc6}
shows that
\begin{equation}
\label{bc9}
\lt|\wt{O}\lt(e_i\rt)-\xi_i\rt|\leq C\ep.
\end{equation}
We define $l_O:\R^n\to \R^n$ to be the linear map such that $l_O\lt(e_i\rt):=\xi_i$ for $i=1,2,\dots n$. Note
$O := Dl_O\in O\lt(n\rt)$. Also, by (\ref{bc9}) we have $ |\wt{O}-O |\leq C\ep$.
In particular $\lt|\zeta_i -O(z_i)\rt| =  |\wt{O}(z_i)-O(z_i) | \le  C \ep$, so that we have
proved \eqref{e10}.

\em Step 5. \rm
Finally, suppose that $\zeta\in \R^n$ and $z\in B_1$ satisfy
$\lt|\lt|\zeta_i -\zeta \rt|-\lt|z_i-z\rt|\rt|\leq \ep$
for all $i \in\lt\{0,1,\dots n\rt\}$.
Then using \eqref{e10} and the fact $O\in O(n)$, we find that
\begin{align*}
\lt|\lt|z_i - l_O^{-1}(\zeta) \rt|-\lt|z_i-z\rt|\rt|
\ &\le \
\lt|\lt|l_O^{-1}(\zeta_i -\zeta) \rt|-\lt|z_i-z\rt|\rt|
+ |l_O^{-1}(\zeta_i) - z_i|\\
&\leq C\ep\
\end{align*}
for all $i \in\lt\{0,1,\dots n\rt\}$.
Arguing exactly as in the proof of \eqref{bc4} in Step 1, we deduce from the above that
\[
| z_i\cdot (l_O^{-1}(\zeta) - z )| \le C \ep
\]
for every $i$.
And  this implies that $|\zeta - l_O(z)| = |l_O^{-1}(\zeta)-z| \le C\ep$; this is proved in Lemma \ref{AX1}
in the next subsection.
\end{proof}

\subsection{Coarea formula into $S^{n-1}$ and bounding the diameter of a simplex}

\begin{a1}
\label{AX3} Let $\Theta_x:\R^n\to S^{n-1}$ be defined by
$\Theta_x\lt(z\rt)=\frac{z-x}{\lt|z-x\rt|}$. Then for any function
$h:\R^n\to \R$ such $H(z) := h(z ) \ |x- z|^{1-n}$ is integrable,

\begin{equation}
\label{xx105} \int_{\psi\in S^{n-1}}\int_{\Theta_x^{-1}\lt(\psi\rt)}
h\lt(z\rt)\, dH^1 z\  d H^{n-1} \psi = \int_{\R^{n}} \frac {h\lt(z\rt)}{|x-z|^{n-1} } dL^n z.
\end{equation}
\end{a1}

\begin{proof}
By a change of variables and Fubini's Theorem,
\begin{align*}
\int_{\R^{n}} \frac {h\lt(z\rt)}{|x-z|^{n-1} } dz &=
\int_0^\infty \int_{z\in \partial B_s(x)} \frac{h\lt(z\rt)}{|x-z|^{n-1} } dH^{n-1}z \, ds\\
&=
\int_0^\infty \int_{\psi\in S^{n-1}} h\lt(s\psi +x \rt)  \, dH^{n-1}\psi \, ds\\
&= \int_{\psi\in S^{n-1}}\int_{\Theta_x^{-1}\lt(\psi\rt)} h\lt(z\rt)
\, dH^1 z \ d H^{n-1} \psi.
\end{align*}
\end{proof}

\begin{a1}
\label{AX1} Let $z_0,z_1,\dots z_n$ be vectors with the property
that $B_b\subset \mathrm{conv}\lt(z_0,z_1,\dots z_n\rt)$,
and  let $S:=\lt\{x:x\cdot z_i\leq 1\text{ for }i=0,1,\dots n\rt\}$.
Then $S\subset B_{n/b}$.
\end{a1}

\begin{proof}[ Proof of Lemma \ref{AX1}]
Fix any $x_0\in S$. Since $b \frac{x_0}{\lt|x_0\rt|}\in \bar
B_b\subset \mathrm{conv}\lt(z_0,z_1,\dots z_n\rt)$ there
exists $\lm_0,\lm_1,\dots \lm_n\in\lt[0,1\rt]$ with $\sum_{i=0}^n
\lm_i z_i=b \frac{x_0}{\lt|x_0\rt|}$. So there must exist $i_0\in
\lt\{0,1,\dots n\rt\}$ such that $z_{i_0}\cdot
\frac{x_0}{\lt|x_0\rt|}\geq \frac{b}{n}$. However as $x_0\in S$ we
have $x_0\cdot z_{i_0}\leq 1$ this gives $\lt|x_0\rt|\leq
\frac{n}{b}$.
\end{proof}

\end{document}